\documentclass[journal]{IEEEtran}
\usepackage[figuresright]{rotating}
\usepackage{amssymb}
\usepackage{amsthm}
\usepackage{multicol}
\usepackage{subfigure}
\usepackage{graphicx}
\usepackage{epstopdf}
\usepackage{fullpage}
\usepackage{latexsym,amsmath}
\usepackage{stmaryrd}
\usepackage{algorithm,algorithmic}
\usepackage{subfigure}
\usepackage{amsfonts}
\usepackage{xcolor,multirow}
\usepackage{cite}
\usepackage{bm}
\usepackage{url}
\usepackage{diagbox}
\usepackage{array}
\usepackage{booktabs}

\newtheorem{lemma}{Lemma}
\newtheorem{theorem}{Theorem}

\ifCLASSINFOpdf
\else
\fi
\begin{document}
	
\title{Quaternion matrix regression for color face recognition}

\author{Jifei~Miao and Kit~Ian~Kou

\thanks{The authors are with the Department of Mathematics, Faculty of Science
	and Technology, University of Macau, Macau 999078, China  (e-mail: jifmiao@163.com;
	kikou@umac.mo)}}


\maketitle
\begin{abstract}
Regression analysis-based approaches have been widely studied for face recognition (FR) in the past several years. More recently, to better deal with some difficult conditions such as occlusions and illumination, nuclear norm based matrix regression methods have been proposed to characterize the low-rank structure of the error image, which generalize the one-dimensional, pixel-based error model to the two-dimensional structure. These methods, however, are inherently devised for grayscale image based FR and without exploiting the color information which is proved beneficial for FR of color face images. Benefiting from quaternion representation, which is capable of encoding the cross-channel correlation of color images, we propose a novel color FR method by formulating the color FR problem as a nuclear norm based quaternion matrix regression (NQMR). We further develop a more robust model called R-NQMR by using the logarithm of the nuclear norm, instead of the original nuclear norm, which adaptively assigns weights on different singular values, and then extend it to deal with the mixed noise. The proposed models, then, are solved using the effective alternating direction multiplier method (ADMM). Experiments on several public face databases demonstrate the superior performance and efficacy of the proposed approaches for color FR, especially for some difficult conditions (occlusion, illumination and mixed noise) over some state-of-the-art regression analysis-based approaches.
\end{abstract}

\IEEEpeerreviewmaketitle

\section{Introduction}
Face recognition (FR) has received widespread attention in the areas of pattern recognition and computer vision due to the growing demand from real-world applications \cite{DBLP:journals/iajit/AzeemSRM14,DBLP:journals/air/LahasanLS19,DBLP:journals/ijprai/BehamR13}. Among numerous FR methods, regression analysis-based approaches have achieved promising results \cite{DBLP:journals/pami/NaseemTB10,DBLP:conf/iccv/ZhangYF11,DBLP:journals/tip/XieYQTZ17,DBLP:journals/pami/YangLQTZX17}, which generally can be classified into two categories, one-dimensional (vectors) regression model \cite{DBLP:journals/pami/NaseemTB10, DBLP:journals/pami/WrightYGSM09, DBLP:conf/iccv/ZhangYF11, DBLP:journals/tip/ZouKW16} and two-dimensional (matrices) regression model \cite{DBLP:journals/tip/XieYQTZ17,DBLP:journals/pami/YangLQTZX17}.

The previous works mainly focus on one-dimensional representation for grayscale FR, \emph{e.g.},
a linear regression classifier (LRC) was proposed in \cite{DBLP:journals/pami/NaseemTB10},  which aims to find a suitable representation of the testing face image, and classify it by checking which class can get the best representation among all the classes. To improve the efficiency, a sparse
representation based classification (SRC) using $L_{1}$-norm regularization on
the coding coefficients and a collaborative representation classifier (CRC) using the $L_{2}$-norm regularization on the coding coefficients schemes for FR were respectively proposed in \cite{DBLP:journals/pami/WrightYGSM09} and \cite{DBLP:conf/iccv/ZhangYF11}. In \cite{DBLP:conf/iccv/ZhangYF11}, the authors also showed that SRC and CRC would get similar results. Although these methods have been successful, there are two major disadvantages. The first one is that they assume the pixels of error follow Gaussian or Laplacian distribution independently \cite{DBLP:journals/tip/FangXLLW15}, hence, the representation error is usually measured by the $L_{2}$-norm or $L_{1}$-norm. Nevertheless, in real-world FR cases because of occlusion, disguise, or illumination variation, the distribution of representation error is more complicated \cite{DBLP:journals/tip/YangZYZ13}. And the second is that these regression methods all use the one-dimensional pixel-based error model to address FR problem, \emph{i.e.}, they consist in transforming the image matrix into vectors in advance, which, therefore, neglects the whole structural information and relationship of the error image \cite{DBLP:journals/pami/YangLQTZX17}.

To overcome these limitations, recently, the authors in \cite{DBLP:journals/pami/YangLQTZX17} found
that occlusion and illumination changes generally lead to a low-rank or approximately low-rank error image, and thus presented the nuclear norm-based matrix regression (NMR) model for FR, which directly characterizes the whole structure of error image with efficient matrix computation. Moreover, the NMR can alleviate the inherent correlations caused by contiguous noises by using the involved singular value decomposition (SVD) \cite{8672214}. The experimental results illustrated that the NMR is robust to real occlusion and illumination. In addition, the authors in \cite{DBLP:journals/tip/XieYQTZ17} substituted the nuclear norm with non-convex function (weighted nuclear norm), and built a robust nuclear norm based matrix regression model (RMR) which can further improve the robustness and deal with mixed noise in FR.

However, the aforementioned methods are originally designed for grayscale FR and fail to
exploit the color information of color face images. Color image contains much more information than does a greyscale image of the same size, which is beneficial to FR \cite{DBLP:conf/icip/TorresRL99,DBLP:journals/mlc/BaoSHYW19,DBLP:journals/tnn/XiaoZ19}. To process color images, one may apply the previous mentioned methods on color channels of the color images independently and separately or concatenate the representations of different color channels into large matrices or form a weighted sum of the three color channel pixels and assign it to corresponding location in new grayscale matrix and, thus, fail to consider the structural correlation information among the color channels \cite{DBLP:journals/tip/LanZT16,DBLP:journals/tip/XuYXZN15}. Nevertheless, such structural correlation information is important for FR \cite{DBLP:journalsper2002}. To preserve the cross-channel correlation for regression analysis-based color FR, the authors in \cite{DBLP:journals/tip/ZouKW16} use quaternion vectors to represent vectorized color images, then extend the traditional CRC and SRC approaches to the quaternion field, \emph{i.e.}, QCRC and QSRC. In addition, the authors in \cite{DBLP:journals/access/ZouKDZT19} developed a quaternion linear representation classification (QLRC) algorithm and a quaternion collaborative representation optimized classifier (QCROC), which are the generalizations of the methods LRC \cite{DBLP:journals/pami/NaseemTB10} and CROC \cite{DBLP:journals/pami/ChiP14} in the quaternion field. The results proved that the quaternion-based algorithms provide a significant improvement over their traditional counterparts.

The existing quaternion-based regression methods mentioned above all use the one-dimensional pixel-based error model, in which, thus, the two major disadvantages described in the previous still exist. In this paper, we first propose a nuclear norm-based quaternion matrix regression (NQMR) model for color FR, which can not only alleviate the two major disadvantages existing in the one-dimensional pixel-based error model but also make the most of color information. In addition, the nuclear norm, \emph{i.e.}, the sum of all singular values, generally over-emphasizes large singular value in practice, which may lead to identification bias \cite{DBLP:journals/tip/XieYQTZ17,DBLP:journals/tnn/ZhengLYBT19}. Furthermore, the image data may contain outliers and noise during image acquisition and transmission procedures \cite{DBLP:journals/jstsp/ChenWZ18}. To address these issues, we further propose a robust nuclear norm-based quaternion matrix regression (R-NQMR) model. Our contributions are listed as follows.
\begin{itemize}
	\item Utilizing quaternion matrices to represent color images and the advantages of the nuclear norm-based matrix regression model, we propose a nuclear norm-based quaternion matrix regression (NQMR) model. We develop an alternating minimization algorithm to iteratively calculate the solution of NQMR in the equivalent real space.
	The optimization, in each iteration, is designed under the framework of the ADMM \cite{DBLP:journals/ftml/BoydPCPE11} which can support the convergence of the algorithm.
	\item To further improve the robustness and practicality, we propose a  novel robust nuclear norm-based quaternion matrix regression (R-NQMR) model. Specifically, we use the logarithm of the nuclear norm to replace the original nuclear norm, which can adaptively assign weights on different singular values \cite{Shuhang17}, and thus be able to effectively approximate the matrix rank. Moreover, we consider mixed noise in our model, adopting the $L_{2}$-norm to describe the Gaussian noise and the $L_{1}$-norm to describe the outliers and non-Gaussian noise, such as salt$\&$pepper noise.
	
	\item The effectiveness of NQMR and R-NQMR is verified by the application of color FR.  Experiments show that NQMR and R-NQMR  can well recognize color face images under difficult conditions such as occlusion and illumination. And R-NQMR is more robust to mixed noise.
\end{itemize}

The remainder of this paper is organized as follows. Section \ref{sec2} provides some notations and preliminaries for quaternion algebra. Section \ref{sec3} first introduces the NMR model. Then, the NQMR model, the R-NQMR model and their optimization algorithms are proposed. Moreover, the computational complexity of the proposed algorithms is also analyzed in this section. The classification strategies for the proposed methods are introduced in Section \ref{sec4}. Section \ref{sec5} provides some experiments to illustrate the performance of our approaches, and compare them with some state-of-the-art related methods. Finally, some conclusions are drawn in Section \ref{sec6}.

\section{Notations and preliminaries}
\label{sec2}
In this section, we first summarize some main notations and then introduce some basic knowledge of quaternion algebra.
\subsection{Notations and definitions}
In this paper, $\mathbb{R}$, and $\mathbb{H}$ respectively denote the set of real numbers and the set of quaternions. A scalar, a vector and a matrix are written as $a$, $\mathbf{a}$ and  $\mathbf{A}$, respectively. A dot (above the variable) is used to denote a quaternion variable (\emph{e.g.,} \cite{DBLP:journals/tip/ZouKW16, DBLP:journals/tip/XuYXZN15}), $\dot{a}$, $\dot{\mathbf{a}}$ and $\dot{\mathbf{A}}$ respectively represent a quaternion scalar, a quaternion vector and a quaternion matrix. $(\cdot)^{\ast}$, $(\cdot)^{T}$, $(\cdot)^{H}$ and $(\cdot)^{-1}$ respectively denote the conjugation, transpose, conjugate transpose and inverse.  $|\cdot|$, $\|\cdot\|_{2}$, $\|\cdot\|_{L_{1}}$, $\|\cdot\|_{F}$, $\|\cdot\|_{2,1}$ and $\|\cdot\|_{\ast}$ are respectively the moduli, the vector $L_{2}$ norm, the matrix $L_{1}$ norm, the Frobenius ( matrix $L_{2}$) norm, the matrix $L_{2,1}$ norm and the nuclear norm.  ${\rm{vec}}(\cdot)$ denotes an operator converting a matrix to a vector by column. And $\mathbf{I}$ denotes the identity matrix with appropriate size.
\subsection{Basic knowledge of quaternion algebra}
Quaternions were discovered in 1843 by W.R. Hamilton \cite{doi:10.1080/14786444408644923}\footnote{Here we just give some fundamental algebraic operations used in our work briefly, which follow the definition in \cite{DBLP:journals/pr/ShaoSWCC14,Girard2007Quaternions}. Readers can find more details on quaternion algebra in the references.}. A quaternion $\dot{q}\in\mathbb{H}$  is a four-dimensional (4D) hypercomplex number and has a Cartesian form given by:
\begin{equation}
\label{equ2}
\dot{q}=q_{0}+q_{1}i+q_{2}j+q_{3}k,
\end{equation}
where $q_{l}\in\mathbb{R}\: (l=0,1,2,3)$ are called its components, and $i, j, k$ are
square roots of -1 and are related through the famous relations:
\begin{align}
\left\{
\begin{array}{lc}
i^{2}=j^{2}=k^{2}=ijk=-1,\\
ij=-ji=k,
jk=-kj=i, 
ki=-ik=j.
\end{array}
\right.
\end{align}
A quaternion $\dot{q}\in\mathbb{H}$ can be decomposed into a real part R$(\dot{q})$ and an imaginary part I$(\dot{q})$ such that $\dot{q}={\rm{R}}(\dot{q})+{\rm{I}}(\dot{q})$, where ${\rm{R}}(\dot{q})=q_{0}$, ${\rm{I}}(\dot{q})=q_{1}i+q_{2}j+q_{3}k$. Then, $\dot{q}\in\mathbb{H}$ will be called a pure quaternion if its real part is null, \emph{i.e.}, if ${\rm{R}}(\dot{q})=0$.
Given two quaternions $\dot{p}$ and $\dot{q}\in\mathbb{H}$, the sum and multiplication of them are respectively:
\begin{equation*}
\label{equ4}
\dot{p}+\dot{q}=(p_{0}+q_{0})+(p_{1}+q_{1})i+(p_{2}+q_{2})j+(p_{3}+q_{3})k
\end{equation*}
and
\begin{align*}
\label{equ5}
\dot{p}\dot{q}=&(p_{0}q_{0}-p_{1}q_{1}-p_{2}q_{2}-p_{3}q_{3})\\
&+(p_{0}q_{1}+p_{1}q_{0}+p_{2}q_{3}-p_{3}q_{2})i\\
&+(p_{0}q_{2}-p_{1}q_{3}+p_{2}q_{0}+p_{3}q_{1})j\\
&+(p_{0}q_{3}+p_{1}q_{2}-p_{2}q_{1}+p_{3}q_{0})k.
\end{align*}
It is noticeable that the multiplication of two quaternions is not
commutative so that in general $\dot{p}\dot{q}\neq\dot{q}\dot{p}$. The conjugate and the modulus of a quaternion $\dot{q}\in\mathbb{H}$ are, respectively, defined as follows:
\begin{align*}
\dot{q}^{\ast}&=q_{0}-q_{1}i-q_{2}j-q_{3}k,\\
|\dot{q}|&=\sqrt{\dot{q}\dot{q}^{\ast}}=\sqrt{q_{0}^{2}+q_{1}^{2}+q_{2}^{2}+q_{3}^{2}}.
\end{align*}
\subsubsection{Quaternion vector and quaternion matrix}
A quaternion vector $\dot{\mathbf{q}}=(\dot{q}_{n})\in\mathbb{H}^{N}$, where $n=1,2,\dots,N$ is a position index, is written as $\dot{\mathbf{q}}=\mathbf{q}_{0}+\mathbf{q}_{1}i+\mathbf{q}_{2}j+\mathbf{q}_{3}k$, where  $\mathbf{q}_{l}\in\mathbb{R}^{N}\: (l=0,1,2,3)$. The $\dot{\mathbf{q}}$ is named a pure quaternion vector when $\mathbf{q}_{0}=\mathbf{0}$. The $L_{2}$-norm of $\dot{\mathbf{q}}$ is defined as $\|\dot{\mathbf{q}}\|_{2}=(\sum_{m=1}^{N}|\dot{q}_{n}|^{2})^{\frac{1}{2}}$. Analogously, a quaternion matrix $\dot{\mathbf{Q}}=(\dot{q}_{mn})\in\mathbb{H}^{M\times N}$, where $m=1,2,\dots,M$ and $n=1,2,\dots,N$ are respectively the row and column indices, is written as $\dot{\mathbf{Q}}=\mathbf{Q}_{0}+\mathbf{Q}_{1}i+\mathbf{Q}_{2}j+\mathbf{Q}_{3}k$, where $\mathbf{Q}_{l}\in\mathbb{R}^{M\times N}\: (l=0,1,2,3)$. The $\dot{\mathbf{Q}}$ is named a pure quaternion matrix when $\mathbf{Q}_{0}=\mathbf{0}$. The  $L_{1}$-norm, $L_{2}$-norm (Frobenius norm) and nuclear norm of $\dot{\mathbf{Q}}$ are respectively defined as $\|\dot{\mathbf{Q}}\|_{L_{1}}=\sum_{m=1}^{M}\sum_{n=1}^{N}|\dot{q}_{mn}|$,  $\|\dot{\mathbf{Q}}\|_{F}=\big(\sum_{m=1}^{M}\sum_{n=1}^{N}|\dot{q}_{mn}|^{2}\big)^{\frac{1}{2}}$ and $\|\dot{\mathbf{Q}}\|_{\ast}=\sum_{k=1}^{r}\sigma_{k}$, where $\sigma_{k}>0$ is the positive singular value of $\dot{\mathbf{Q}}$.
\subsubsection{Two quaternion operators $\mathcal{P}$ and $\mathcal{Q}$}
The noncommutativity of quaternion multiplication makes it quite difficult to handle quaternion optimization problems. Hence, we first introduce two useful quaternion operators, which can convert the problem from quaternion number field to
real number field \cite{DBLP:journals/tip/XuYXZN15,DBLP:journals/tip/ZouKW16}. For any quaternion matrix  $\dot{\mathbf{Q}}=\mathbf{Q}_{0}+\mathbf{Q}_{1}i+\mathbf{Q}_{2}j+\mathbf{Q}_{3}k\in\mathbb{H}^{M\times N}$, the operator $\mathcal{P}$ is defined as $\mathbb{H}^{M\times N}\to\mathbb{R}^{4M\times 4N}$ such that
\begin{equation}\small
\label{qtorm}
\mathcal{P}(\dot{\mathbf{Q}})=\left(\begin{array}{cccc}
	\mathbf{Q}_{0}&  -\mathbf{Q}_{1}&  -\mathbf{Q}_{2}& -\mathbf{Q}_{3} \\ 
	\mathbf{Q}_{1}&  \mathbf{Q}_{0}&  -\mathbf{Q}_{3} & \mathbf{Q}_{2} \\ 
    \mathbf{Q}_{2}&  \mathbf{Q}_{3}&   \mathbf{Q}_{0}& -\mathbf{Q}_{1}  \\ 
	\mathbf{Q}_{3}&  -\mathbf{Q}_{2}&  \mathbf{Q}_{1}& \mathbf{Q}_{0}
\end{array} \right).
\end{equation}
For any quaternion vector  $\dot{\mathbf{q}}=\mathbf{q}_{0}+\mathbf{q}_{1}i+\mathbf{q}_{2}j+\mathbf{q}_{3}k\in\mathbb{H}^{N}$, the operator $\mathcal{Q}$ is defined as $\mathbb{H}^{N}\to\mathbb{R}^{4N}$ such that
\begin{equation}
\label{qtorv}
\mathcal{Q}(\dot{\mathbf{q}})=(\mathbf{q}_{0}^{T},\mathbf{q}_{1}^{T},\mathbf{q}_{2}^{T},\mathbf{q}_{3}^{T})^{T}.
\end{equation}
$\mathcal{P}(\dot{\mathbf{Q}})$ and $\mathcal{Q}(\dot{\mathbf{q}})$  are uniquely determined by $\dot{\mathbf{Q}}$ and $\dot{\mathbf{q}}$, respectively. The inverse operators of $\mathcal{P}$ and $\mathcal{Q}$ are respectively denoted as $\mathcal{P}^{-1}$ and $\mathcal{Q}^{-1}$.
There are some properties of $\mathcal{P}$ and $\mathcal{Q}$ (\emph{see} Theorem \ref{theorem1}).
\begin{theorem}
	\label{theorem1}
	Let $\dot{\mathbf{P}}$ and $\dot{\mathbf{Q}}\in\mathbb{H}^{M\times N}$, $\dot{\mathbf{p}}\in\mathbb{H}^{N}$. Then,
	\begin{enumerate}
		\item [1)] $\mathcal{P}$ and $\mathcal{Q}$ are linear.
		\item [2)] $\|\mathcal{Q}(\dot{\mathbf{p}})\|_{2}=\|\dot{\mathbf{p}}\|_{2}$.
		\item [3)] $\mathcal{Q}(\dot{\mathbf{Q}}\dot{\mathbf{p}})=\mathcal{P}(\dot{\mathbf{Q}})\mathcal{Q}(\dot{\mathbf{p}})$.
		\item [4)] $\|\dot{\mathbf{Q}}\dot{\mathbf{p}}\|_{2}=\|\mathcal{P}(\dot{\mathbf{Q}})\mathcal{Q}(\dot{\mathbf{p}})\|_{2}$.
		\item [5)] $\|\mathcal{P}(\dot{\mathbf{Q}})\|_{F}=4\|\dot{\mathbf{Q}}\|_{F}$.
		\item [6)] $\|\mathcal{P}(\dot{\mathbf{Q}})\|_{\ast}=4\|\dot{\mathbf{Q}}\|_{\ast}$.
	\end{enumerate}
\end{theorem}
The proof of the Theorem \ref{theorem1} can be found in APPENDIX \ref{appendices1}. From the properties $5)$ and $6)$, we can see that the minimization problems of $\mathop{\rm{min}}\limits_{\dot{\mathbf{Q}}}\|\dot{\mathbf{Q}}\|_{F}$ and $\mathop{\rm{min}}\limits_{\dot{\mathbf{Q}}}\|\dot{\mathbf{Q}}\|_{\ast}$ are respectively equivalent to the minimization problems of $\mathop{\rm{min}}\limits_{\dot{\mathbf{Q}}}\|\mathcal{P}(\dot{\mathbf{Q}})\|_{F}$ and $\mathop{\rm{min}}\limits_{\dot{\mathbf{Q}}}\|\mathcal{P}(\dot{\mathbf{Q}})\|_{\ast}$. 

\section{Methodology}
\label{sec3}
In this section, we first introduce the nuclear norm-based matrix regression (NMR) model. Then, the nuclear norm-based quaternion matrix regression (NQMR) model, the robust nuclear norm-based quaternion matrix regression (R-NQMR) model and their optimization algorithms are proposed.
\subsection{Nuclear norm-based matrix regression}
Given a set of $L$ training image matrices $\mathbf{A}_{1}, \mathbf{A}_{2},\dots, \mathbf{A}_{L}\in\mathbb{R}^{M\times N}$ and a query image matrix $\mathbf{B}\in\mathbb{R}^{M\times N}$ which can be represented linearly as
\begin{equation}
\label{equ6}
\mathbf{B}=x_{1}\mathbf{A}_{1}+x_{2}\mathbf{A}_{2}+\ldots+x_{L}\mathbf{A}_{L}+\mathbf{E},
\end{equation}
where $\mathbf{x}=(x_{1},x_{2},\ldots,x_{L})\in\mathbb{R}^{L}$ is a set of representation coefficients, and $\mathbf{E}\in\mathbb{R}^{M\times N}$ is error image. For convenience, we generally define $\mathbf{A}(\mathbf{x})=x_{1}\mathbf{A}_{1}+x_{2}\mathbf{A}_{2}+\ldots+x_{L}\mathbf{A}_{L}$. Then the formula (\ref{equ6}) can be rewritten as
\begin{equation}
\label{equ7}
\mathbf{B}=\mathbf{A}(\mathbf{x})+\mathbf{E}.
\end{equation}
Based on the assumption that the error image usually has low rank \cite{DBLP:journals/pami/YangLQTZX17}, the representation coefficients are estimated by solving
the following nuclear norm approximation problem
\begin{equation}
\mathop{\rm{min}}\limits_{\mathbf{x}}\|\mathbf{A}(\mathbf{x})-B\|_{\ast}.
\end{equation}
Furthermore, to avoid overfitting, by adding an $L_{2}$-norm regularization term on representation coefficients, the NMR model can be expressed as follows:
\begin{equation}
\mathop{\rm{min}}\limits_{\mathbf{x}}\|\mathbf{A}(\mathbf{x})-\mathbf{B}\|_{\ast}+\frac{\lambda}{2}\|\mathbf{x}\|_{2}^{2},
\end{equation}
where $\lambda$ is a positive parameter.
Then, the regularized matrix regression problem can be solved by the ADMM algorithm \cite{DBLP:journals/pami/YangLQTZX17,DBLP:journals/ijon/DengLWD18}.

\subsection{Nuclear norm-based quaternion matrix regression}
To process color FR problem, we model each RGB color image as a pure quaternion matrix $\dot{\mathbf{P}}\in\mathbb{H}^{M\times N}$ written as
\begin{equation*}
\dot{\mathbf{P}}=\mathbf{0}+\mathbf{P}_{r}i+\mathbf{P}_{g}j+\mathbf{P}_{b}k,
\end{equation*}
where $\mathbf{P}_{r}$, $\mathbf{P}_{g}$ and $\mathbf{P}_{b}$ are, respectively, the red, green and blue channel of $\dot{\mathbf{P}}$. 

Similar to NMR, given a set of $L$ training color image matrices $\dot{\mathbf{A}}_{1}, \dot{\mathbf{A}}_{2},\dots, \dot{\mathbf{A}}_{L}\in\mathbb{H}^{M\times N}$ and a query color image matrix $\dot{\mathbf{B}}\in\mathbb{H}^{M\times N}$ which can be represented linearly as
\begin{equation}
\label{equ8}
\dot{\mathbf{B}}=\dot{x}_{1}\dot{\mathbf{A}}_{1}+\dot{x}_{2}\dot{\mathbf{A}}_{2}+\ldots+\dot{x}_{L}\dot{\mathbf{A}}_{L}+\dot{\mathbf{E}},
\end{equation}
where $\dot{x}_{1},\dot{x}_{2},\ldots,\dot{x}_{L}\in\mathbb{H}$ is a set of representation coefficients, and $\dot{\mathbf{E}}\in\mathbb{H}^{M\times N}$ is error image. The objective of NQMR is to estimate the representation coefficients by solving the following regularized nuclear norm approximation problem
\begin{equation}
\label{equ9}
\mathop{\rm{min}}\limits_{\dot{\mathbf{x}}}\|\dot{\mathbf{A}}(\dot{\mathbf{x}})-\dot{\mathbf{B}}\|_{\ast}+\frac{\lambda}{2}\|\dot{\mathbf{x}}\|_{2}^{2},
\end{equation}
where $\dot{\mathbf{A}}(\dot{\mathbf{x}})=\dot{x}_{1}\dot{\mathbf{A}}_{1}+\dot{x}_{2}\dot{\mathbf{A}}_{2}+\ldots+\dot{x}_{L}\dot{\mathbf{A}}_{L}$ and $\lambda$ is a positive parameter.

We adopt the ADMM  method to solve the problem (\ref{equ9}) in this paper, which can guarantee the convergence of the algorithm \cite{DBLP:journals/corr/LinCM10}. We first convert (\ref{equ9}) to the following equivalent problem:
\begin{equation}
\label{equ10}
\left\{
\begin{array}{c}
\mathop{\rm{min}}\limits_{\dot{\mathbf{x}}}\|\dot{\mathbf{E}}\|_{\ast}+\frac{\lambda}{2}\|\dot{\mathbf{x}}\|_{2}^{2},\\
\text{s.t.} \quad \dot{\mathbf{A}}(\dot{\mathbf{x}})-\dot{\mathbf{B}}= \dot{\mathbf{E}}.
\end{array}
\right.
\end{equation}
Then, the augmented Lagrange function is given by
\begin{align}
\label{equ11}
\mathcal{L}_{\mu}(\dot{\mathbf{E}},\dot{\mathbf{x}}, \dot{\mathbf{\Lambda}})=&\|\dot{\mathbf{E}}\|_{\ast}+\frac{\lambda}{2}\|\dot{\mathbf{x}}\|_{2}^{2}+\langle\dot{\mathbf{\Lambda}},\dot{\mathbf{A}}(\dot{\mathbf{x}})-\dot{\mathbf{B}}-\dot{\mathbf{E}}\rangle \nonumber\\
&+\frac{\mu}{2}\|\dot{\mathbf{A}}(\dot{\mathbf{x}})-\dot{\mathbf{B}}-\dot{\mathbf{E}}\|_{F}^{2},
\end{align}
where $\mu>0$ is a penalty parameter, $\dot{\mathbf{\Lambda}}$ is the Lagrange multiplier.

To update each variable in (\ref{equ11}), in the $\tau+1$th iteration, we perform the following steps:
\begin{itemize}
	\item \textbf{Step 1:} Given $\dot{\mathbf{E}}^{\tau}$ and $\dot{\mathbf{\Lambda}}^{\tau}$, updating $\dot{\mathbf{x}}$ by
	\begin{equation*}
	\dot{\mathbf{x}}^{\tau+1}=\mathop{{\rm{arg\, min}}}\limits_{\dot{\mathbf{x}}}\:\mathcal{L}_{\mu}(\dot{\mathbf{E}}^{\tau},\dot{\mathbf{x}}, \dot{\mathbf{\Lambda}}^{\tau}).
	\end{equation*}
	\item \textbf{Step 2:} Given $\dot{\mathbf{x}}^{\tau+1}$ and $\dot{\mathbf{\Lambda}}^{\tau}$, updating $\dot{\mathbf{E}}$ by
	\begin{equation*}
	\dot{\mathbf{E}}^{\tau+1}=\mathop{{\rm{arg\, min}}}\limits_{\dot{\mathbf{E}}}\:\mathcal{L}_{\mu}(\dot{\mathbf{E}},\dot{\mathbf{x}}^{\tau+1}, \dot{\mathbf{\Lambda}}^{\tau}).
	\end{equation*}
	\item \textbf{Step 3:}  Given $\dot{\mathbf{x}}^{\tau+1}$ and $\dot{\mathbf{E}}^{\tau+1}$, updating $\dot{\mathbf{\Lambda}}$ by
	\begin{equation*}
	\dot{\mathbf{\Lambda}}^{\tau+1}=\dot{\mathbf{\Lambda}}^{\tau}+\mu(\dot{\mathbf{A}}(\dot{\mathbf{x}}^{\tau+1})-\dot{\mathbf{B}}-\dot{\mathbf{E}}^{\tau+1}).
	\end{equation*}	
\end{itemize}
\textbf{Step 3} is a proximal minimization step of the Lagrange multiplies $\dot{\mathbf{\Lambda}}$, \textbf{Step 1} and \textbf{Step 2} are key steps to solve the optimization problem.

For \textbf{Step 1}, $\dot{\mathbf{x}}^{(\tau+1)}$ is the optimal solution of the following problem:
\begin{align}
\label{equ12}
\dot{\mathbf{x}}^{\tau+1}=&\mathop{{\rm{arg\, min}}}\limits_{\dot{\mathbf{x}}}\frac{\lambda}{2}\|\dot{\mathbf{x}}\|_{2}^{2}+\langle\dot{\mathbf{\Lambda}}^{\tau},\dot{\mathbf{A}}(\dot{\mathbf{x}})-\dot{\mathbf{B}}-\dot{\mathbf{E}}^{\tau}\rangle \nonumber\\
&+\frac{\mu}{2}\|\dot{\mathbf{\Lambda}}^{\tau},\dot{\mathbf{A}}(\dot{\mathbf{x}})-\dot{\mathbf{B}}-\dot{\mathbf{E}}^{\tau}\|_{F}^{2} \nonumber\\
=&\mathop{{\rm{arg\, min}}}\limits_{\dot{\mathbf{x}}}\frac{\mu}{2}\|\dot{\mathbf{A}}(\dot{\mathbf{x}})-(\dot{\mathbf{B}}+\dot{\mathbf{E}}^{\tau}-\frac{1}{\mu}\dot{\mathbf{\Lambda}}^{\tau})\|_{F}^{2} \nonumber\\
&+\frac{\lambda}{2}\|\dot{\mathbf{x}}\|_{2}^{2}.
\end{align}
Denote $\dot{\mathbf{H}}=({\rm{vec}}(\dot{\mathbf{A}}_{1}), {\rm{vec}}(\dot{\mathbf{A}}_{2}),\ldots, {\rm{vec}}(\dot{\mathbf{A}}_{L}))$ and $\dot{\mathbf{g}}^{\tau}= {\rm{vec}}(\dot{\mathbf{B}}+\dot{\mathbf{E}}^{\tau}-\frac{1}{\mu}\dot{\mathbf{\Lambda}}^{\tau})$. Then (\ref{equ12}) is equivalent to
\begin{equation}
\label{equ13}
\dot{\mathbf{x}}^{\tau+1}=\mathop{{\rm{arg\, min}}}\limits_{\dot{\mathbf{x}}}\left( \| \dot{\mathbf{H}}\dot{\mathbf{x}}-\dot{\mathbf{g}}^{\tau}\|_{2}^{2}+\frac{\lambda}{\mu}\|\dot{\mathbf{x}}\|_{2}^{2}\right).
\end{equation}
By using the properties of operators $\mathcal{P}$ and $\mathcal{Q}$ in Theorem \ref{theorem1}, we transform the model (\ref{equ13}) into the following equivalent form
\begin{align}
\label{equ14}
\dot{\mathbf{x}}^{\tau+1}=&\mathop{{\rm{arg\, min}}}\limits_{\dot{\mathbf{x}}}\left(\| \mathcal{P}(\dot{\mathbf{H}})\mathcal{Q}(\dot{\mathbf{x}})-\mathcal{Q}(\dot{\mathbf{g}}^{\tau})\|_{2}^{2}\right.\nonumber\\
&\left.+\frac{\lambda}{\mu}\|\mathcal{Q}(\dot{\mathbf{x}})\|_{2}^{2}\right),
\end{align}
which is a standard ridge regression model, and has the following closed-form solution \cite{DBLP:journals/tip/ZouKW16}:
\begin{equation}\small
\label{equ15}
\dot{\mathbf{x}}^{\tau+1}\!=\mathcal{Q}^{-1}\!\left(\!\!\bigg(\mathcal{P}(\dot{\mathbf{H}})^{T}\mathcal{P}(\dot{\mathbf{H}})+\frac{\lambda}{\mu}\mathbf{I}\bigg)^{-1}\!\!\mathcal{P}(\dot{\mathbf{H}})^{T}\mathcal{Q}(\dot{\mathbf{g}}^{\tau}) \right).
\end{equation}

For \textbf{Step 2}, $\dot{\mathbf{E}}^{(\tau+1)}$ is the optimal solution of the following problem:
\begin{align}
\label{equ16}
\dot{\mathbf{E}}^{\tau+1}=&\mathop{{\rm{arg\, min}}}\limits_{\dot{\mathbf{E}}}\|\dot{\mathbf{E}}\|_{\ast}+\langle\dot{\mathbf{\Lambda}}^{\tau},\dot{\mathbf{A}}(\dot{\mathbf{x}}^{\tau+1})-\dot{\mathbf{B}}-\dot{\mathbf{E}}\rangle \nonumber\\
&+\frac{\mu}{2}\|\dot{\mathbf{\Lambda}}^{\tau},\dot{\mathbf{A}}(\dot{\mathbf{x}}^{\tau+1})-\dot{\mathbf{B}}-\dot{\mathbf{E}}\|_{F}^{2}\nonumber\\
=&\mathop{{\rm{arg\, min}}}\limits_{\dot{\mathbf{E}}}\frac{1}{2}\|\dot{\mathbf{E}}-\big(\dot{\mathbf{A}}(\dot{\mathbf{x}}^{\tau+1})-\dot{\mathbf{B}}+\frac{1}{\mu}\dot{\mathbf{\Lambda}}^{\tau}\big) \|_{F}^{2}\nonumber\\
&+\frac{1}{\mu}\|\dot{\mathbf{E}}\|_{\ast}.
\end{align}
Based on the properties of operator $\mathcal{P}$ in Theorem \ref{theorem1}, (\ref{equ16}) is equivalent to the following problem:
\begin{align}
\label{equ17}
\dot{\mathbf{E}}^{\tau+1}\!&=\mathcal{P}^{-1}\left(\mathcal{P}(\dot{\mathbf{E}}^{(\tau+1)})\right)\nonumber\\
&=\mathop{{\rm{arg\, min}}}\limits_{\dot{\mathbf{E}}}\frac{1}{2}\|\!\mathcal{P}(\dot{\mathbf{E}})\!-\!\mathcal{P}\big(\dot{\mathbf{A}}(\!\dot{\mathbf{x}}^{\tau+1}\!)\!-\!\dot{\mathbf{B}}\!+\!\frac{1}{\mu}\dot{\mathbf{\Lambda}}^{\tau}\big)\!\|_{F}^{2}\nonumber\\
&\quad +\frac{1}{\mu}\|\mathcal{P}(\dot{\mathbf{E}})\|_{\ast}\nonumber\\
&=\mathcal{P}^{-1}\left(D_{\frac{1}{\mu}}\bigg(\mathcal{P}\big(\dot{\mathbf{A}}(\dot{\mathbf{x}}^{\tau+1})-\dot{\mathbf{B}}+\frac{1}{\mu}\dot{\mathbf{\Lambda}}^{\tau}\big)\bigg)\right),
\end{align}
where $D_{\gamma}(\mathbf{M})$ is the matrix singular value thresholding operator \cite{DBLP:journals/siamjo/CaiCS10}: $D_{\gamma}(\mathbf{M}):=\mathbf{U}{\rm{diag}}(\bar{\sigma}_{i},i=1,\ldots,r)\mathbf{V}^{T}$, $\mathbf{M}=\mathbf{U}{\rm{diag}}(\sigma_{i},i=1,\ldots,r)\mathbf{V}^{T}$ is the SVD of matrix $\mathbf{M}$ and $\bar{\sigma}_{i}={\rm{max}}\{\sigma_{i}-\gamma,0\}$. 

\subsection{Robust nuclear norm-based quaternion matrix regression}
In real-world applications, images may be corrupted by noise during the image acquisition and transmission procedures. Thus, we tend to consider the noise terms in our model. In addition, we use the logarithm of the nuclear norm to replace the original nuclear norm (\emph{e.g.}, \cite{DBLP:journals/tcyb/XieLFL19}), which can adaptively assigns weights on different singular values \cite{Shuhang17}, and thus can more effectively approximate the matrix rank. Hence, we use the following more robust model to replace (\ref{equ10}):
\begin{equation}
\label{equ19}
\left\{
\begin{array}{lc}
\mathop{\rm{min}}\limits_{\dot{\mathbf{x}}}\omega L(\dot{\mathbf{E}}_{0},\epsilon)+\alpha\|\dot{\mathbf{E}}_{1}\|_{L_{1}}+\beta\|\dot{\mathbf{E}}_{2}\|_{F}^{2}+\eta\|\dot{\mathbf{x}}\|_{2}^{2},\\
\text{s.t.} \quad \dot{\mathbf{A}}(\dot{\mathbf{x}})-\dot{\mathbf{B}}= \dot{\mathbf{E}}_{0}+\dot{\mathbf{E}}_{1}+\dot{\mathbf{E}}_{2},
\end{array}
\right.
\end{equation}
where $\epsilon$ is a small positive number, $L(\dot{\mathbf{E}}_{0},\epsilon)=\sum_{i}{\rm{log}}(\sigma_{i}(\dot{\mathbf{E}}_{0})+\epsilon)$, $\sigma_{i}(\dot{\mathbf{E}}_{0})$ represents the $i$th singular value of quaternion  matrix $\dot{\mathbf{E}}_{0}$, $\omega$, $\alpha$, $\beta$ and $\eta$ are positive numbers. The $\dot{\mathbf{E}}_{0}$ is the low-rank quaternion  matrix term, which is the clean color image data.  The $\dot{\mathbf{E}}_{1}$ is the sparse error term denoting outliers and non-Gaussian noise, The $\dot{\mathbf{E}}_{2}$ denotes the Gaussian noise. Using the ADMM framework, the augmented Lagrange function of (\ref{equ19}) is given by
\begin{align}
\label{equ20}
&\mathcal{F}_{\mu}(\dot{\mathbf{E}}_{0},\dot{\mathbf{E}}_{1},\dot{\mathbf{E}}_{2},\dot{\mathbf{x}}, \dot{\mathbf{\Lambda}}) \nonumber\\
&=\omega L(\dot{\mathbf{E}}_{0},\epsilon)+\alpha\|\dot{\mathbf{E}}_{1}\|_{L_{1}}+\beta\|\dot{\mathbf{E}}_{2}\|_{F}^{2} +\eta\|\dot{\mathbf{x}}\|_{2}^{2}\nonumber\\
&\quad +\langle\dot{\mathbf{\Lambda}},\dot{\mathbf{A}}(\dot{\mathbf{x}})+\dot{\mathbf{E}}_{0}+\dot{\mathbf{E}}_{1}+\dot{\mathbf{E}}_{2}-\dot{\mathbf{B}}\rangle \nonumber\\
&\quad +\frac{\mu}{2}\|\dot{\mathbf{A}}(\dot{\mathbf{x}})+\dot{\mathbf{E}}_{0}+\dot{\mathbf{E}}_{1}+\dot{\mathbf{E}}_{2}-\dot{\mathbf{B}}\|_{F}^{2},
\end{align}
where $\mu>0$ is a penalty parameter, $\dot{\mathbf{\Lambda}}$ is the Lagrange multiplier. 

To update each variable in (\ref{equ20}), we perform the following iterative scheme:

\subsubsection{Update $\dot{\mathbf{x}}$} Fixing other variables except for $\dot{\mathbf{x}}$ in (\ref{equ20}), we solve the following problem:
\begin{align}
\label{equ34}
\dot{\mathbf{x}}^{\tau+1}=&\mathop{{\rm{arg\, min}}}\limits_{\dot{\mathbf{x}}}\eta\|\dot{\mathbf{x}}\|_{2}^{2}\nonumber\\
&+\langle\dot{\mathbf{\Lambda}}^{\tau},\dot{\mathbf{A}}(\dot{\mathbf{x}})+\dot{\mathbf{E}}_{0}^{\tau}+\dot{\mathbf{E}}_{1}^{\tau}+\dot{\mathbf{E}}_{2}^{\tau}-\dot{\mathbf{B}}\rangle \nonumber\\
&+\frac{\mu}{2}\|\dot{\mathbf{A}}(\dot{\mathbf{x}})+\dot{\mathbf{E}}_{0}^{\tau}+\dot{\mathbf{E}}_{1}^{\tau}+\dot{\mathbf{E}}_{2}^{\tau}-\dot{\mathbf{B}}\|_{F}^{2}\nonumber\\
=&\mathop{{\rm{arg\, min}}}\limits_{\dot{\mathbf{x}}}\frac{\eta}{\mu}\|\dot{\mathbf{x}}\|_{2}^{2}\nonumber\\
&+\frac{1}{2}\|\dot{\mathbf{A}}(\dot{\mathbf{x}})-\big(\dot{\mathbf{B}}-\dot{\mathbf{E}}_{0}^{\tau}-\dot{\mathbf{E}}_{1}^{\tau}\nonumber\\
&-\dot{\mathbf{E}}_{2}^{\tau}-\frac{1}{\mu}\mathbf{\Lambda}^{\tau}\big)\|_{F}^{2}.
\end{align}
Denote $\dot{\mathbf{y}}^{\tau}= {\rm{vec}}(\dot{\mathbf{B}}-\dot{\mathbf{E}}_{0}^{\tau}-\dot{\mathbf{E}}_{1}^{\tau}-\dot{\mathbf{E}}_{2}^{\tau}-\frac{1}{\mu}\mathbf{\Lambda}^{\tau})$. Then (\ref{equ34}) is equivalent to
\begin{equation*}
\label{equ35}
\dot{\mathbf{x}}^{\tau+1}=\mathop{{\rm{arg\, min}}}\limits_{\dot{\mathbf{x}}}\frac{\eta}{\mu}\|\dot{\mathbf{x}}\|_{2}^{2}+\frac{1}{2}\|\dot{\mathbf{H}}\dot{\mathbf{x}}-\dot{\mathbf{y}}^{\tau}\|_{2}^{2},
\end{equation*}
which is equivalent to
\begin{align}
\label{equ36}
\dot{\mathbf{x}}^{\tau+1}=&\mathop{{\rm{arg\, min}}}\limits_{\dot{\mathbf{x}}}\left(\| \mathcal{P}(\dot{\mathbf{H}})\mathcal{Q}(\dot{\mathbf{x}})-\mathcal{Q}(\dot{\mathbf{y}}^{\tau})\|_{2}^{2}\right.\nonumber\\
&\left.+\frac{\eta}{\mu}\|\mathcal{Q}(\dot{\mathbf{x}})\|_{2}^{2}\right),
\end{align}
and has the following closed-form solution:
\begin{equation}\small
\label{equ37}
\dot{\mathbf{x}}^{\tau+1}\!=\mathcal{Q}^{-1}\!\left(\!\!\bigg(\mathcal{P}(\dot{\mathbf{H}})^{T}\mathcal{P}(\dot{\mathbf{H}})+\frac{\eta}{\mu}\mathbf{I}\bigg)^{-1}\!\!\mathcal{P}(\dot{\mathbf{H}})^{T}\mathcal{Q}(\dot{\mathbf{y}}^{\tau}) \right).
\end{equation}

\subsubsection{Update $\dot{\mathbf{E}}_{0}$} Fixing other variables except for $\dot{\mathbf{E}}_{0}$ in (\ref{equ20}), we solve the following problem:
\begin{align}
\label{equ21}
\dot{\mathbf{E}}_{0}^{\tau+1}=&\mathop{{\rm{arg\, min}}}\limits_{\dot{\mathbf{E}}_{0}}
\omega L(\dot{\mathbf{E}}_{0},\epsilon)\nonumber\\
&+\langle\dot{\mathbf{\Lambda}}^{\tau},\dot{\mathbf{A}}(\dot{\mathbf{x}}^{\tau+1})+\dot{\mathbf{E}}_{0}+\dot{\mathbf{E}}_{1}^{\tau}+\dot{\mathbf{E}}_{2}^{\tau}-\dot{\mathbf{B}}\rangle \nonumber\\
&+\frac{\mu}{2}\|\dot{\mathbf{A}}(\dot{\mathbf{x}}^{\tau+1})+\dot{\mathbf{E}}_{0}+\dot{\mathbf{E}}_{1}^{\tau}+\dot{\mathbf{E}}_{2}^{\tau}-\dot{\mathbf{B}}\|_{F}^{2}.
\end{align}
Note that $L(\dot{\mathbf{E}}_{0},\epsilon)=\sum_{i}{\rm{log}}(\sigma_{i}(\dot{\mathbf{E}}_{0})+\epsilon)$, and $\sigma_{i}(\dot{\mathbf{E}}_{0})\in \mathbb{R}$, thus ${\rm{log}}(\sigma_{i}(\dot{\mathbf{E}}_{0})+\epsilon)$ can be approximated by its first-order Taylor expansion: ${\rm{log}}(\sigma_{i}(\dot{\mathbf{E}}_{0})+\epsilon)\approx{\rm{log}}(\sigma_{i}(\dot{\mathbf{E}}_{0}^{\tau})+\epsilon)+\langle \nabla{\rm{log}}(\sigma_{i}(\dot{\mathbf{E}}_{0}^{\tau})+\epsilon), \sigma_{i}(\dot{\mathbf{E}}_{0})-\sigma_{i}(\dot{\mathbf{E}}_{0}^{\tau}) \rangle$, where $\sigma_{i}(\dot{\mathbf{E}}_{0}^{\tau})$ is the value obtained in the $\tau$th iteration. We can see that
${\rm{log}}(\sigma_{i}(\dot{\mathbf{E}}_{0}^{\tau})+\epsilon)$ and $\sigma_{i}(\dot{\mathbf{E}}_{0}^{\tau})/(\sigma_{i}(\dot{\mathbf{E}}_{0}^{\tau})+\epsilon)$
are constants and can be ignored since they do not affect the minimization problem, \emph{i.e.}, the minimization of ${\rm{log}}(\sigma_{i}(\dot{\mathbf{E}}_{0})+\epsilon)$ can be approximated by the minimization of $\sigma_{i}(\dot{\mathbf{E}}_{0})/(\sigma_{i}(\dot{\mathbf{E}}_{0}^{\tau})+\epsilon)$. As a result, the problem $\min
L(\dot{\mathbf{E}}_{0},\epsilon)$  can be replaced by $\min\sum_{i}(\sigma_{i}(\dot{\mathbf{E}}_{0})/(\sigma_{i}(\dot{\mathbf{E}}_{0}^{\tau})+\epsilon))$. Then, we rewrite the problem (\ref{equ21}) as
\begin{align}
\label{equ22}
\dot{\mathbf{E}}_{0}^{\tau+1}
=&\mathop{{\rm{arg\, min}}}\limits_{\dot{\mathbf{E}}_{0}}
\frac{\omega}{\mu}\sum_{i}\frac{\sigma_{i}(\dot{\mathbf{E}}_{0})}{\sigma_{i}(\dot{\mathbf{E}}_{0}^{\tau})+\epsilon}\nonumber\\
&+\!\frac{1}{2}\|\dot{\mathbf{E}}_{0}\!-\!\big(\dot{\mathbf{B}}\!-\dot{\mathbf{A}}(\dot{\mathbf{x}}^{\tau+1})\!\!-\!\dot{\mathbf{E}}_{1}^{\tau}\!-\!\dot{\mathbf{E}}_{2}^{\tau}\!-\!\frac{1}{\mu}\mathbf{\Lambda}^{\tau}\big)\!\|_{F}^{2}.
\end{align}
We can observe that $\frac{\omega}{\mu}\sum_{i}\frac{\sigma_{i}(\dot{\mathbf{E}}_{0})}{\sigma_{i}(\dot{\mathbf{E}}_{0}^{\tau})+\epsilon}$ is actually the weighted nuclear norm of $\dot{\mathbf{E}}_{0}$ with adaptive weight $\varGamma_{i}^{\tau}:=\frac{\omega/\mu}{\sigma_{i}(\dot{\mathbf{E}}_{0}^{\tau})+\epsilon}$ of each singular value $\sigma_{i}(\dot{\mathbf{E}}_{0})$. Thus, for convenience, we denote $\frac{\omega}{\mu}\sum_{i}\frac{\sigma_{i}(\dot{\mathbf{E}}_{0})}{\sigma_{i}(\dot{\mathbf{E}}_{0}^{\tau})+\epsilon}$ as $\|\dot{\mathbf{E}}_{0}\|_{(\varGamma^{\tau},\ast)}$. Then, (\ref{equ22}) is equivalent to
\begin{align}
\label{equ23}
\dot{\mathbf{E}}_{0}^{\tau+1}=&\mathop{{\rm{arg\, min}}}\limits_{\dot{\mathbf{E}}_{0}}
\|\dot{\mathbf{E}}_{0}\|_{(\varGamma^{\tau},\ast)}\nonumber\\
&+\frac{1}{2}\|\mathcal{P}(\dot{\mathbf{E}}_{0})-\mathcal{P}\big(\dot{\mathbf{B}}-\dot{\mathbf{A}}(\dot{\mathbf{x}}^{\tau+1})-\dot{\mathbf{E}}_{1}^{\tau}\nonumber\\
&-\dot{\mathbf{E}}_{2}^{\tau}-\frac{1}{\mu}\mathbf{\Lambda}^{\tau}\big)\|_{F}^{2}.
\end{align}
The problem (\ref{equ23}) has the following closed-form solution \cite{DBLP:conf/cvpr/GuZZF14}:
\begin{align}
\label{equ24}
\dot{\mathbf{E}}_{0}^{\tau+1}=&\mathcal{P}^{-1}\left(\Psi_{\varGamma^{\tau}}\big(\mathcal{P}\big(\dot{\mathbf{B}}-\dot{\mathbf{A}}(\dot{\mathbf{x}}^{\tau+1})-\dot{\mathbf{E}}_{1}^{\tau}\right.\nonumber\\
&\left.-\dot{\mathbf{E}}_{2}^{\tau}-\frac{1}{\mu}\mathbf{\Lambda}^{\tau}\big)\big)\right),
\end{align}
where $\Psi_{\mathbf{s}}(\mathbf{C}):=\mathbf{U}{\rm{diag}}(\widetilde{\sigma}_{1},\ldots,\widetilde{\sigma}_{r})\mathbf{V}^{T}$, $\mathbf{C}=\mathbf{U}{\rm{diag}}(\sigma_{1},\ldots,\sigma_{r})\mathbf{V}^{T}$ is the SVD of matrix $\mathbf{C}$ and
\begin{align*}
\widetilde{\sigma}_{i}=\left\{
\begin{array}{lc}
\sigma_{i}-s_{i} &\text{if} \ \sigma_{i}-s_{i}>0,\\
0 &\text{otherwise}.
\end{array}
\right.
\end{align*}
\subsubsection{Update $\dot{\mathbf{E}}_{1}$} Fixing other variables except for $\dot{\mathbf{E}}_{1}$ in (\ref{equ20}), we solve the following problem:
\begin{align}
\label{equ26}
\dot{\mathbf{E}}_{1}^{\tau+1}=&\mathop{{\rm{arg\, min}}}\limits_{\dot{\mathbf{E}}_{1}}\alpha\|\dot{\mathbf{E}}_{1}\|_{L_{1}}\nonumber\\
&+\langle\dot{\mathbf{\Lambda}}^{\tau},\dot{\mathbf{A}}(\dot{\mathbf{x}^{\tau+1}})+\dot{\mathbf{E}}_{0}^{\tau+1}+\dot{\mathbf{E}}_{1}+\dot{\mathbf{E}}_{2}^{\tau}-\dot{\mathbf{B}}\rangle \nonumber\\
&+\frac{\mu}{2}\|\dot{\mathbf{A}}(\dot{\mathbf{x}^{\tau+1}})+\dot{\mathbf{E}}_{0}^{\tau+1}+\dot{\mathbf{E}}_{1}+\dot{\mathbf{E}}_{2}^{\tau}-\dot{\mathbf{B}}\|_{F}^{2}\nonumber\\
=&\mathop{{\rm{arg\, min}}}\limits_{\dot{\mathbf{E}}_{1}}\frac{\alpha}{\mu}\|\dot{\mathbf{E}}_{1}\|_{L_{1}}\nonumber\\
&+\frac{1}{2}\|\dot{\mathbf{E}}_{1}-\big(\dot{\mathbf{B}}-\dot{\mathbf{A}}(\dot{\mathbf{x}}^{\tau+1})-\dot{\mathbf{E}}_{0}^{\tau+1}-\dot{\mathbf{E}}_{2}^{\tau}\nonumber\\
&-\frac{1}{\mu}\mathbf{\Lambda}^{\tau}\big)\|_{F}^{2}.
\end{align}
Before going further, we first denote a useful operator $\mathcal{R}:\mathbb{H}^{M\times N}\to \mathbb{R}^{4\times MN}$ such that
\begin{equation}
\label{equ27}
\mathcal{R}\big(\dot{\mathbf{Q}})=\big({\rm{vec}}(\mathbf{Q}_{0}),{\rm{vec}}(\mathbf{Q}_{1}),{\rm{vec}}(\mathbf{Q}_{2}),{\rm{vec}}(\mathbf{Q}_{3})\big)^{T}.   
\end{equation}
There are some properties of $\mathcal{R}$ (\emph{see} Theorem \ref{theorem2}).
\begin{theorem}
	\label{theorem2}
	Let $\dot{\mathbf{P}}$ and $\dot{\mathbf{Q}}\in\mathbb{H}^{M\times N}$. Then,
	\begin{enumerate}
		\item [1)] $\mathcal{R}$ is linear, \emph{i.e.}, $\mathcal{R}(\dot{\mathbf{P}}+\dot{\mathbf{Q}})=\mathcal{R}(\dot{\mathbf{P}})+\mathcal{R}(\dot{\mathbf{Q}})$.
		\item [2)] $\|\mathcal{R}(\dot{\mathbf{Q}})\|_{F}=\|\dot{\mathbf{Q}}\|_{F}$.
		\item [3)] For any quaternion matrix $\dot{\mathbf{Q}}$, we have
		\begin{equation}
		\label{equ28}
		\|\dot{\mathbf{Q}}\|_{L_{1}}=\|\mathcal{R}\big(\dot{\mathbf{Q}})\|_{2,1},
		\end{equation}
		where $\|\mathbf{M}\|_{2,1}$ denotes the $L_{2,1}$ norm of matrix $\mathbf{M}$, which is defined as $\|\mathbf{M}\|_{2,1}:=\sum_{j}\|\mathbf{M}(:,j)\|_{2}$. 
	\end{enumerate}
\end{theorem}
According to the definition of $\mathcal{R}$ and with some real matrix calculation, the properties $1)-3)$ are obviously true.
Thus, based on the Theorem \ref{theorem2}, problem (\ref{equ26}) is equivalent to
\begin{align}
\label{equ29}
\dot{\mathbf{E}}_{1}^{\tau+1}=&\mathop{{\rm{arg\, min}}}\limits_{\dot{\mathbf{E}}_{1}}
\frac{\alpha}{\mu}\|\mathcal{R}(\dot{\mathbf{E}}_{1})\|_{2,1}\nonumber\\
&+\frac{1}{2}\|\mathcal{R}(\dot{\mathbf{E}}_{1})-\mathcal{R}\big(\dot{\mathbf{B}}-\dot{\mathbf{A}}(\dot{\mathbf{x}}^{\tau+1})-\dot{\mathbf{E}}_{0}^{\tau+1}\nonumber\\
&-\dot{\mathbf{E}}_{2}^{\tau}-\frac{1}{\mu}\mathbf{\Lambda}^{\tau}\big)\|_{F}^{2}.
\end{align}
$\mathcal{R}(\dot{\mathbf{E}}_{1}^{\tau+1})$ in (\ref{equ29}) can be found based on the following Lemma \ref{lemma1}
\begin{lemma}\cite{DBLP:journals/tnn/XiaoZ19}
\label{lemma1}
The optimal solution $\hat{\mathbf{Z}}$ to the problem
\begin{equation*}
\hat{\mathbf{Z}}=\mathop{{\rm{arg\, min}}}\limits_{\mathbf{Z}}\lambda\|\mathbf{Z}\|_{2,1}+\frac{1}{2}\|\mathbf{Z}-\mathbf{M}\|_{F}^{2},
\end{equation*}
satisfies 
\begin{align*}
\hat{\mathbf{Z}}(:,j)=\left\{
\begin{array}{lc}
\frac{\|\mathbf{M}(:,j)\|_{2}-\frac{\alpha}{\mu}}{\|\mathbf{M}(:,j)\|_{2}}\mathbf{M}(:,j), &\|\mathbf{M}(:,j)\|_{2}>\frac{\alpha}{\mu},\\
\mathbf{0}, &\text{otherwise}.
\end{array}
\right.
\end{align*}
\end{lemma}
Then, 
\begin{equation}
\label{equ30}
\dot{\mathbf{E}}_{1}^{\tau+1}=\mathcal{R}^{-1}\big(\mathcal{R}(\dot{\mathbf{E}}_{1}^{\tau+1})\big). 
\end{equation}
\subsubsection{Update $\dot{\mathbf{E}}_{2}$} Fixing other variables except for $\dot{\mathbf{E}}_{2}$ in (\ref{equ20}), we solve the following problem:
\begin{align}
\label{equ31}
\dot{\mathbf{E}}_{2}^{\tau+1}=&\mathop{{\rm{arg\, min}}}\limits_{\dot{\mathbf{E}}_{2}}\beta\|\dot{\mathbf{E}}_{2}\|_{F}^{2}\nonumber\\
&+\langle\dot{\mathbf{\Lambda}}^{\tau},\dot{\mathbf{A}}(\dot{\mathbf{x}^{\tau+1}})+\dot{\mathbf{E}}_{0}^{\tau+1}+\dot{\mathbf{E}}_{1}^{\tau+1}+\dot{\mathbf{E}}_{2}-\dot{\mathbf{B}}\rangle \nonumber\\
&+\frac{\mu}{2}\|\dot{\mathbf{A}}(\dot{\mathbf{x}^{\tau+1}})+\dot{\mathbf{E}}_{0}^{\tau+1}+\dot{\mathbf{E}}_{1}^{\tau+1}+\dot{\mathbf{E}}_{2}-\dot{\mathbf{B}}\|_{F}^{2}\nonumber\\
=&\mathop{{\rm{arg\, min}}}\limits_{\dot{\mathbf{E}}_{2}}\frac{\beta}{\mu}\|\dot{\mathbf{E}}_{2}\|_{F}^{2}\nonumber\\
&+\frac{1}{2}\|\dot{\mathbf{E}}_{2}-\big(\dot{\mathbf{B}}-\dot{\mathbf{A}}(\dot{\mathbf{x}}^{\tau+1})-\dot{\mathbf{E}}_{0}^{\tau+1}-\dot{\mathbf{E}}_{1}^{\tau+1}\nonumber\\
&-\frac{1}{\mu}\mathbf{\Lambda}^{\tau}\big)\|_{F}^{2},
\end{align}
which is equivalent to
\begin{align}
\label{equ32}
\dot{\mathbf{E}}_{2}^{\tau+1}=&\mathop{{\rm{arg\, min}}}\limits_{\dot{\mathbf{E}}_{2}}\frac{\beta}{\mu}\|\mathcal{P}(\dot{\mathbf{E}}_{2})\|_{F}^{2}\nonumber\\
&+\frac{1}{2}\|\mathcal{P}(\dot{\mathbf{E}}_{2})-\mathcal{P}\big(\dot{\mathbf{B}}-\dot{\mathbf{A}}(\dot{\mathbf{x}}^{\tau+1})\nonumber\\
&-\dot{\mathbf{E}}_{0}^{\tau+1}-\dot{\mathbf{E}}_{1}^{\tau+1}-\frac{1}{\mu}\mathbf{\Lambda}^{\tau}\big)\|_{F}^{2}.
\end{align}
This is a standard least squares regression problem with closed-form solution:
\begin{align}
\label{equ33}
\dot{\mathbf{E}}_{2}^{\tau+1}=&\mathcal{P}^{-1}\bigg((\beta+\mu)^{-1}\mu\mathcal{P}\big(\dot{\mathbf{B}}-\dot{\mathbf{A}}(\dot{\mathbf{x}}^{\tau+1})
-\dot{\mathbf{E}}_{0}^{\tau+1}\nonumber\\
&-\dot{\mathbf{E}}_{1}^{\tau+1}-\frac{1}{\mu}\mathbf{\Lambda}^{\tau}\big)\bigg).
\end{align}

Finally, $\dot{\mathbf{\Lambda}}$ is updated as
\begin{align}
\label{equ38}
\dot{\mathbf{\Lambda}}^{\tau+1}=&\dot{\mathbf{\Lambda}}^{\tau}+\mu\left(\dot{\mathbf{A}}(\dot{\mathbf{x}}^{\tau+1})+\dot{\mathbf{E}}_{0}^{\tau+1}+\dot{\mathbf{E}}_{1}^{\tau+1}\right.\nonumber\\
&\left.+\dot{\mathbf{E}}_{2}^{\tau+1}-\dot{\mathbf{B}}\right).
\end{align}

The whole algorithms for NQMR and R-NQMR are respectively summarized in Table \ref{tab_algorithm1} and Table \ref{tab_algorithm2}.
\begin{table}[htbp]
	\caption{Nuclear norm based quaternion matrix regression (\textbf{NQMR}) algorithm.}
	\hrule
	\label{tab_algorithm1}
	\begin{algorithmic}[1]
		\REQUIRE A set of $L$ training color image matrices $\dot{\mathbf{A}}_{1}, \dot{\mathbf{A}}_{2},\dots, \dot{\mathbf{A}}_{L}\in\mathbb{H}^{M\times N}$ and a query color image matrix $\dot{\mathbf{B}}\in\mathbb{H}^{M\times N}$, the parameters $\lambda$, $\mu$ and $\varepsilon_{rel}$, \emph{Max-iter}.
		\STATE \textbf{Initialize}  $\dot{\mathbf{E}}^{0}=\dot{\mathbf{0}}$ and $\dot{\mathbf{\Lambda}}^{0}=\dot{\mathbf{0}}$.
		\STATE Let $\dot{\mathbf{H}}=({\rm{vec}}(\dot{\mathbf{A}}_{1}), {\rm{vec}}(\dot{\mathbf{A}}_{2}),\ldots, {\rm{vec}}(\dot{\mathbf{A}}_{L}))$ and compute $\mathbf{Z}:=\big(\mathcal{P}(\dot{\mathbf{H}})^{T}\mathcal{P}(\dot{\mathbf{H}})+\frac{\lambda}{\mu}\mathbf{I}\big)^{-1}\mathcal{P}(\dot{\mathbf{H}})^{T}$.
		\STATE \textbf{Repeat}
		\STATE Update $\dot{\mathbf{x}}: \dot{\mathbf{x}}^{\tau+1}=\mathcal{Q}^{-1}\left(\mathbf{Z}\mathcal{Q}(\dot{\mathbf{g}}^{\tau})\right)$.
		\STATE Update $\dot{\mathbf{E}}$ using equation (\ref{equ17}).
		\STATE Update $\dot{\mathbf{\Lambda}}: \dot{\mathbf{\Lambda}}^{\tau+1}=\dot{\mathbf{\Lambda}}^{\tau}+\mu(\dot{\mathbf{A}}(\dot{\mathbf{x}}^{\tau+1})-\dot{\mathbf{B}}-\dot{\mathbf{E}}^{\tau+1})$.
		\STATE \textbf{Until} $|Dual_{(\tau)}-Dual_{(\tau+1)}|<\varepsilon_{rel}$ or the number of iterations$>\emph{Max-iter}$, where $Dual_{(\tau)}=\|\dot{\mathbf{A}}(\dot{\mathbf{x}}^{\tau+1})-\dot{\mathbf{B}}-\dot{\mathbf{E}}^{\tau+1}\|_{F}$.
		\ENSURE $\dot{\mathbf{x}}^{\tau+1}$.
	\end{algorithmic}
	\hrule
\end{table}

\begin{table}[htbp]
	\caption{Robust nuclear norm-based quaternion matrix regression (\textbf{R-NQMR}) algorithm.}
	\hrule
	\label{tab_algorithm2}
		\begin{algorithmic}[1]
		\REQUIRE A set of $L$ training color image matrices $\dot{\mathbf{A}}_{1}, \dot{\mathbf{A}}_{2},\dots, \dot{\mathbf{A}}_{L}\in\mathbb{H}^{M\times N}$ and a query color image matrix $\dot{\mathbf{B}}\in\mathbb{H}^{M\times N}$, the parameters $\omega$, $\alpha$, $\beta$, $\eta$, $\mu$ and $\varepsilon_{rel}$, \emph{Max-iter}.
		\STATE \textbf{Initialize}   $\dot{\mathbf{E}}_{0}^{0}=\dot{\mathbf{0}}$, $\dot{\mathbf{E}}_{1}^{0}=\dot{\mathbf{0}}$, $\dot{\mathbf{E}}_{2}^{0}=\dot{\mathbf{0}}$ and $\dot{\mathbf{\Lambda}}^{0}=\dot{\mathbf{0}}$.
		\STATE Let $\dot{\mathbf{H}}=({\rm{vec}}(\dot{\mathbf{A}}_{1}), {\rm{vec}}(\dot{\mathbf{A}}_{2}),\ldots, {\rm{vec}}(\dot{\mathbf{A}}_{L}))$ and compute $\mathbf{L}:=\big(\mathcal{P}(\dot{\mathbf{H}})^{T}\mathcal{P}(\dot{\mathbf{H}})+\frac{\eta}{\mu}\mathbf{I}\big)^{-1}\mathcal{P}(\dot{\mathbf{H}})^{T}$.
		\STATE \textbf{Repeat}
		\STATE Update $\dot{\mathbf{x}}: \dot{\mathbf{x}}^{\tau+1}=\mathcal{Q}^{-1}\left(\mathbf{L}\mathcal{Q}(\dot{\mathbf{y}}^{\tau})\right)$.
		\STATE Update $\dot{\mathbf{E}}_{0}$ based on equation (\ref{equ24}).
		\STATE Update $\dot{\mathbf{E}}_{1}$ based on the Lemma \ref{lemma1} and equation (\ref{equ30}).
		\STATE Update $\dot{\mathbf{E}}_{2}$ based on equations (\ref{equ32}) and (\ref{equ33}).
		\STATE Update $\dot{\mathbf{\Lambda}}$ based on equation (\ref{equ38}).	
		\STATE \textbf{Until} $|Dual_{(\tau)}-Dual_{(\tau+1)}|<\varepsilon_{rel}$ or the number of iterations$>\emph{Max-iter}$, where $Dual_{(\tau)}=\|\dot{\mathbf{A}}(\dot{\mathbf{x}}^{\tau+1})-\dot{\mathbf{B}}-\dot{\mathbf{E}}_{0}^{\tau+1}
		-\dot{\mathbf{E}}_{1}^{\tau+1}-\dot{\mathbf{E}}_{2}^{\tau+1}\|_{F}$.
		\ENSURE $\dot{\mathbf{x}}^{\tau+1}$.
	\end{algorithmic}
	\hrule
\end{table}
\subsection{The computational complexity}
Suppose that there is a set of $L$ training color image matrices $\dot{\mathbf{A}}_{1}, \dot{\mathbf{A}}_{2},\dots, \dot{\mathbf{A}}_{L}\in\mathbb{H}^{M\times N}$. For convenience,
we assume that $M\leq N$. Note that $\mathbf{Z}:=\big(\mathcal{P}(\dot{\mathbf{H}})^{T}\mathcal{P}(\dot{\mathbf{H}})+\frac{\lambda}{\mu}\mathbf{I}\big)^{-1}\mathcal{P}(\dot{\mathbf{H}})^{T}$ and 
$\mathbf{L}:=\big(\mathcal{P}(\dot{\mathbf{H}})^{T}\mathcal{P}(\dot{\mathbf{H}})+\frac{\eta}{\mu}\mathbf{I}\big)^{-1}\mathcal{P}(\dot{\mathbf{H}})^{T}$  are fixed in each iteration, thus they can be calculated and stored beforehand. NQMR: The cost of updating $\dot{\mathbf{x}}$ is determined by the matrix-vector multiplication $\mathbf{Z}\mathcal{Q}(\dot{\mathbf{g}})$, which is $\mathcal{O}(16LMN)$. The cost of updating $\dot{\mathbf{E}}$  is determined by performing SVD on $4M\times 4N$ real-valued matrix, which is $\mathcal{O}(64MN^{2})$. Therefore, in each iteration, the total computational complexity of NQMR algorithm is about $\mathcal{O}(16LMN+64MN^{2})$. R-NQMR: The cost of updating $\dot{\mathbf{x}}$ and $\dot{\mathbf{E}}_{0}$ is approximately equal to that of updating $\dot{\mathbf{x}}$ and $\dot{\mathbf{E}}$ in NQMR algorithm. The updating of $\dot{\mathbf{E}}_{2}$ only involves the product of matrix and single constant, so the cost can be ignored relative to the whole cost. Hence, the total computational complexity of R-NQMR algorithm is about $\mathcal{O}(16LMN+64MN^{2})$.

\section{Classification  strategies}
\label{sec4}
In this section, we introduce the NQMR and R-NQMR based classification  strategies. Given a set of $L$ training color image matrices $\dot{\mathbf{A}}_{1}, \dot{\mathbf{A}}_{2},\dots, \dot{\mathbf{A}}_{L}\in\mathbb{H}^{M\times N}$ from different classes and a new query color image matrix $\dot{\mathbf{B}}\in\mathbb{H}^{M\times N}$, we use $L$ training color image matrices to represent $\dot{\mathbf{B}}$ and obtain the representation coefficients based on the NQMR or R-NQMR algorithm.

Based on the optimal solution, \emph{e.g.}, $\dot{\mathbf{x}}_{*}$, we reconstruct the color image matrix $\dot{\mathbf{B}}$ as $\hat{\dot{\mathbf{B}}}=\dot{\mathbf{A}}(\dot{\mathbf{x}}_{*})$. Define operator $\delta_{k}: \mathbb{H}^{n}\to \mathbb{H}^{n}$ such that $\delta_{k}(\dot{\mathbf{x}}_{*})$ is a vector whose only nonzero entries are the entries of $\dot{\mathbf{x}}_{*}$ associated with the $k$-th class. Then, we reconstruct the color image matrix $\dot{\mathbf{B}}$ in the class $k$ as $\hat{\dot{\mathbf{B}}}_{k}=\dot{\mathbf{A}}(\delta_{k}(\dot{\mathbf{x}}_{*}))$. For NQMR model, the class-dependent reconstruction error is defined by
\begin{equation*}
\label{equ39}
r_{k}(\dot{\mathbf{B}})=\|\hat{\dot{\mathbf{B}}}-\hat{\dot{\mathbf{B}}}_{k}\|_{\ast}=\|\dot{\mathbf{A}}(\dot{\mathbf{x}}_{*})-\dot{\mathbf{A}}(\delta_{k}(\dot{\mathbf{x}}_{*}))\|_{\ast}.
\end{equation*}
For R-NQMR model, the class-dependent reconstruction error is defined by
\begin{align}
\label{equ40}
r_{k}(\dot{\mathbf{B}})&=\|\hat{\dot{\mathbf{B}}}-\hat{\dot{\mathbf{B}}}_{k}\|_{(\varGamma^{ f inal},\ast)}\nonumber\\
&=\|\dot{\mathbf{A}}(\dot{\mathbf{x}}_{*})-\dot{\mathbf{A}}(\delta_{k}(\dot{\mathbf{x}}_{*}))\|_{(\varGamma^{final},\ast)},
\end{align}
where $\varGamma^{final}$ is the final weighted vector.
 
After computing all the class-dependent reconstruction errors (suppose there are total $K$ classes), the query color image matrix $\dot{\mathbf{B}}$ is assigned to the class with minimal error
\begin{equation}
\label{equ41}
\text{identify}(\dot{\mathbf{B}})=\mathop{{\rm{arg\, min}}}\limits_{k\in K}r_{k}(\dot{\mathbf{B}}).
\end{equation}

\section{Experimental results}
\label{sec5}
In this section, several color face databases are selected to validate the performance of the proposed NQMR and R-NQMR methods. And we compare them with several existing state-of-the-art regression analysis-based approaches including LRC \cite{DBLP:journals/pami/NaseemTB10}, QLRC \cite{DBLP:journals/access/ZouKDZT19}, CRC \cite{DBLP:conf/iccv/ZhangYF11}, QCRC \cite{DBLP:journals/tip/ZouKW16}, CROC \cite{DBLP:journals/pami/ChiP14}, QCROC \cite{DBLP:journals/access/ZouKDZT19}, NMR \cite{DBLP:journals/pami/YangLQTZX17}, RMR and S-RMR ($L_{2}$-norm regularization pose on the coding coefficients) \cite{DBLP:journals/tip/XieYQTZ17}. To represent color face images, QLRC, QCRC and QCROC use quaternion vectors, while NQMR and R-NQMR use quaternion matrices. Other algorithms utilize the real value vectors or real value matrices to represent color images by using weighted sum of the three color channel pixels and assign it to corresponding location in new grayscale matrices (\emph{i.e.}, ``rgb2gray'' function in MATLAB).

\textbf{Parameter Settings:} For NQMR in TABLE  \ref{tab_algorithm1}, $\lambda=1$, $\mu=1$ and $\varepsilon_{rel}=10^{-4}$. For R-NQMR in TABLE \ref{tab_algorithm2}, $\eta=1$, $\mu=1$ and $\varepsilon_{rel}=10^{-4}$. We pick the values for $\omega$, $\alpha$ and $\beta$, say $(0.01, 0.1, 1, 10)$, and we choose these three parameters by the highest recognition rate. All the parameters are determined empirically. The parameter settings of all the other methods follow the author's suggestions. 

\subsection{Experiments on the AR face database \cite{Benavente1998}} 
\begin{figure}[htbp]
	\centering
	\subfigure[]{\includegraphics[width=7.5cm,height=1.7cm]{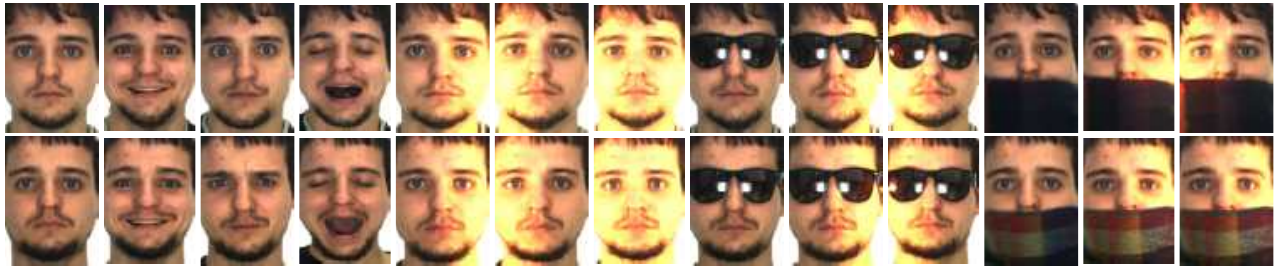}}\\
	\subfigure[]{\includegraphics[width=7.5cm,height=1.7cm]{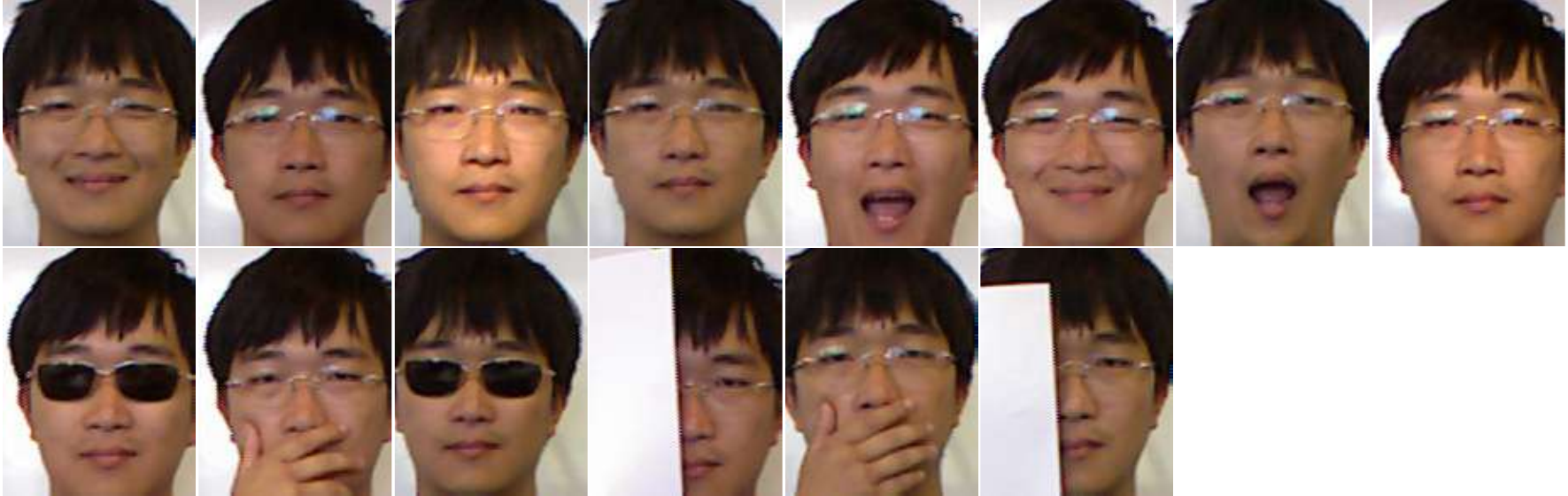}}\\
	\subfigure[]{\includegraphics[width=7.5cm,height=1.7cm]{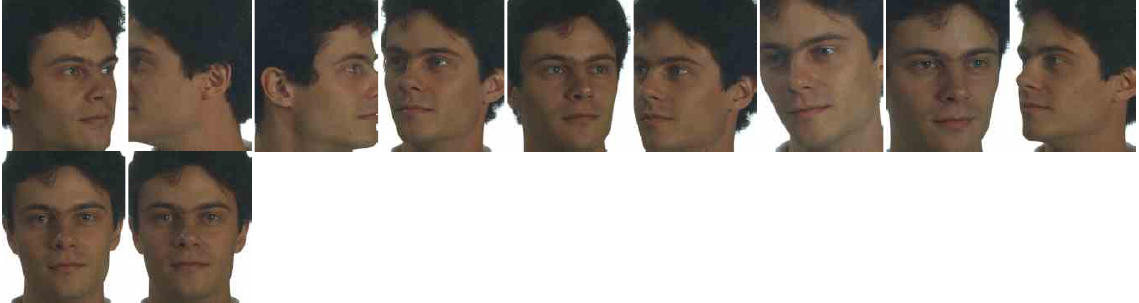}}
	\caption{Sample color face images of one person from (a) AR database, (b) EURECOM kinect database, (c) Color FERET database.}
	\label{fig1}
\end{figure}

The AR face database contains over $4,000$ color images corresponding to $126$ people's faces ($70$ men and $56$ women). Images feature frontal view faces with different facial expressions, illumination conditions, and occlusions (sunglasses and scarf). Each person participated in two sessions, separated by two weeks ($14$ days) time. We select a random subset of AR containing $100$ subjects, \emph{i.e.}, $2600$ images in total are used. The images of one person from AR are shown in Fig. \ref{fig1}(a). All face images are resized to $42\times 30$ pixels.

\subsubsection{Performance on realistic illumination and occlusions} This experiment is divided into two cases. The first case (\textbf{case 1}) chooses the first $8$ images (with various facial expressions) from both sessions of each person (\emph{e.g.}, the first four images of the first and the second rows in Fig. \ref{fig1}(a)) to form the training set.
The test set is formed by $18$ images with various illumination and occlusions (glasses or scarves)
from both sessions of each person (\emph{e.g.}, the $5$th to the $13$th images of the first and the second rows in Fig. \ref{fig1}(a)). The second case (\textbf{case 2}) chooses $8$ images with different illumination from both sessions of each person (\emph{e.g.}, the $4$th to the $7$th images of the first and the second rows in Fig. \ref{fig1}(a)) to form the training set. The test set is formed by $12$ occluded (glasses or scarves) images
from both sessions of each person (\emph{e.g.}, the $8$th to the $13$th images of the first and the second rows in Fig. \ref{fig1}(a)).

\begin{table*}	
\caption{Recognition rates (\%) of different methods on the AR database for the two cases (the best results are highlighted in bold for each case).}
\centering
\begin{tabular}{|c|c|c|c|c|c|c|c|c|c|c|c|}
	\hline 
	\diagbox{Cases}{Rates (\%)}{Methods}&LRC&  QLRC&  CRC&  QCRC&  CROC& QCROC & NMR & RMR & S-RMR & NQMR & R-NQMR \\ 
	\hline 
	case 1&28.17 &31.56  &64.11  &71.89  &46.78  & 52.28 & 77.89 & 84.39 & 84.41 & 79.11 &\textbf{87.77}  \\ 
	\hline 
	case 2&19.08& 22.17 & 64.50 & 74.92 & 39.25 &45.58 & 80.58 & 85.00 & 85.16 & 82.17 & \textbf{89.16} \\ 
	\hline 
\end{tabular} 
\label{tab3}
\end{table*}

From TABLE \ref{tab3}, we find that, for the two cases, matrix regression based approaches (NMR, RMR, S-RMR, NQMR and R-NQMR) generally show better performance than that of vector regression based approaches (LRC, QLRC, CRC, QCRC, CROC and QCROC). The performance of quaternion-based approaches (QLRC, QCRC, QCROC and NQMR) have obvious advantages over their real version counterparts (LRC, CRC, CROC and NMR). And the recognition rate of R-NQMR is much better than that of other methods.

\subsubsection{Performance on random block occlusions} 
\begin{figure}[htbp]
	\centering
	\includegraphics[width=8cm,height=4cm]{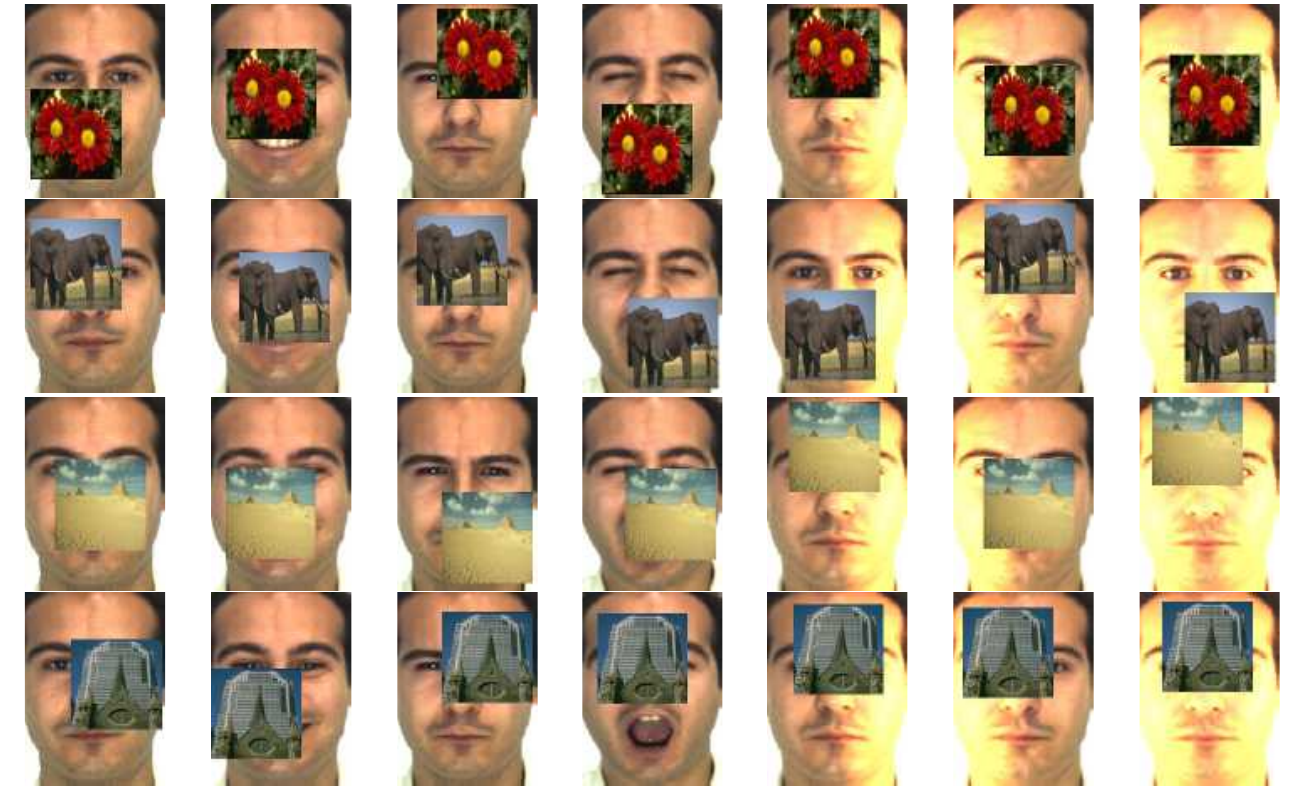}
	\caption{Example color face images of one person from the AR database simultaneously contain different random block occlusions (from top to bottom:  flower, elephant, desert and building).}
	\label{fig2}
\end{figure}

\begin{figure}[htbp]
	\centering
	\includegraphics[width=8.5cm,height=5cm]{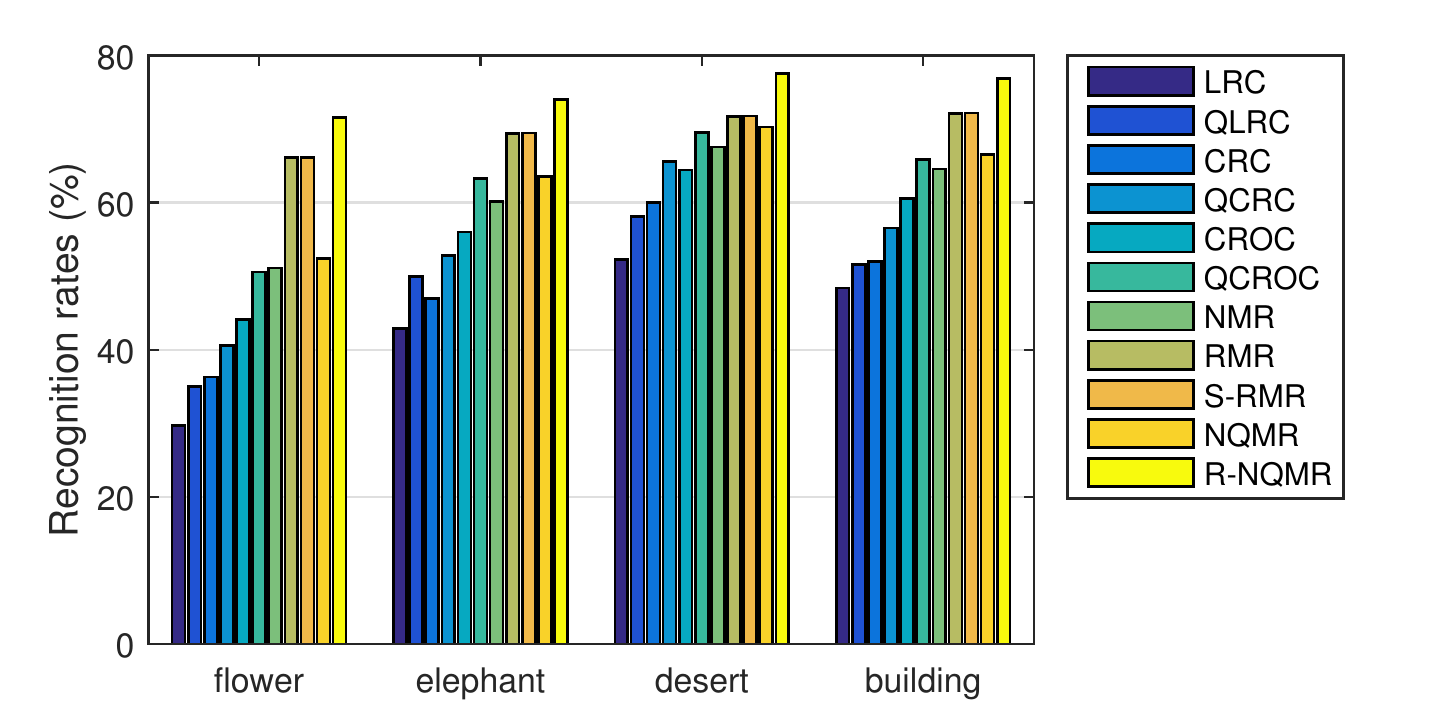}
	\caption{Recognition rates (\%) under various contents of block occlusions (from left to right:  flower, elephant, desert, building).}
	\label{fig3}
\end{figure}

\begin{figure}[htbp]
	\centering
	\includegraphics[width=8cm,height=4cm]{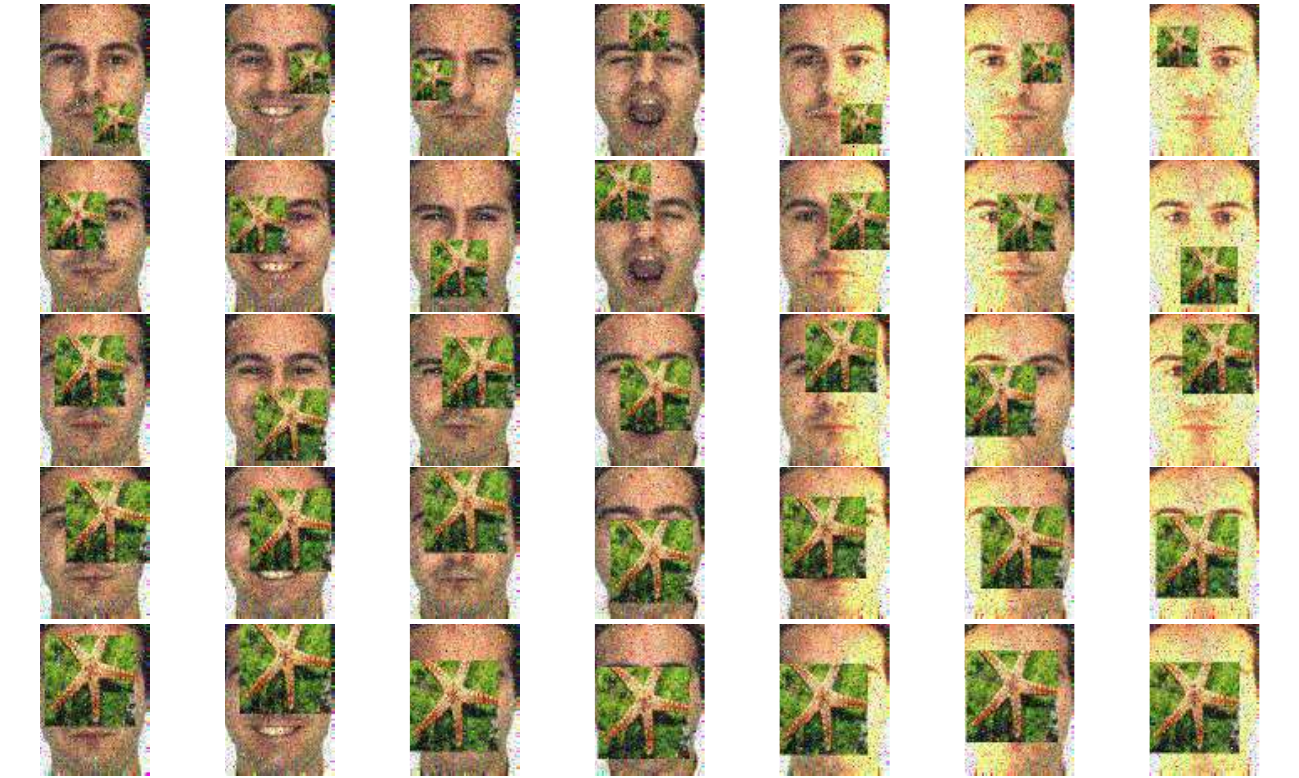}
	\caption{Example color face images of one person from the AR database with different sizes of block occlusions (from top to bottom:  $10\%$, $20\%$, $30\%$, $40\%$, $50\%$) and mixed noise (``slat \& pepper'' noise and Gaussian white noise).}
	\label{fig4}
\end{figure}

\begin{figure*}[htbp]
	\centering
	\includegraphics[width=18.7cm,height=6cm]{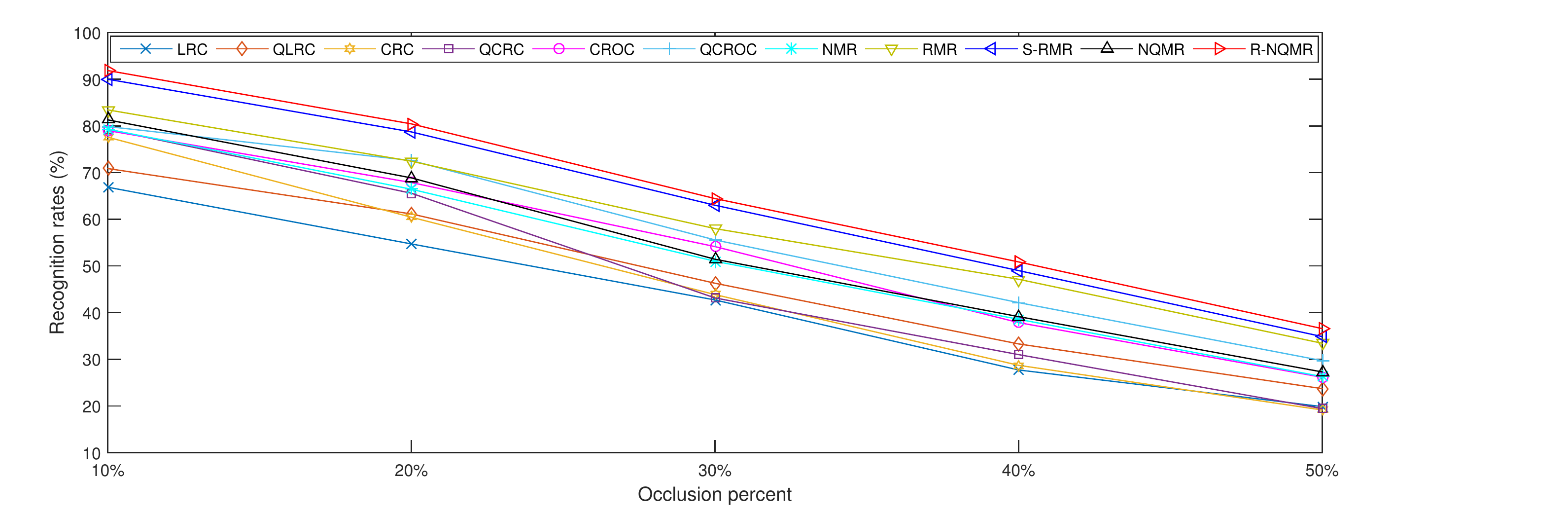}
	\caption{Recognition rates (\%) under varying portions of occlusions (from left to right: $10\%$, $20\%$, $30\%$, $40\%$ and $50\%$).}
	\label{fig5}
\end{figure*}

In this experiment, we add synthetic occlusions to evaluate the performance of competing methods. And the experiment is divided into two parts. In the \textbf{first part}, $7$ nonoccluded color face images from session one  (\emph{e.g.}, the first to the $7$th images of the first row in Fig. \ref{fig1}(a)) are used for training, and $7$ nonoccluded images from session two (\emph{e.g.}, the first to the $7$th images of the second row in Fig. \ref{fig1}(a)) are added with an unrelated color image for testing. Each occlusion is scaled to $30\%$ of the size of the color face images and is imposed at a random block. Shown as in Fig. \ref{fig2}, $4$ different occlusions (flower, elephant, desert and building) are selected. According to Fig. \ref{fig3}, R-NQMR shows the highest recognition rate for all kind of occlusions. In the \textbf{second part}, we use the same training set as in the first part. The testing images are still 
 the $7$ nonoccluded images from session two
but with different sizes (from $10\%$ to $50\%$ of the size of the color face images) of uncorrelated blocks and mixed noise (including ``salt \& pepper'' noise with $0.1$ probabilities and Gaussian white noise with zero mean and $0.01$ variance). Example test images of one person are shown as Fig. \ref{fig4}. From Fig. \ref{fig5}, we can see that R-NQMR has consistently the best performance
among all percents of occlusions.

\subsection{Experiments on the EURECOM kinect database \cite{IEEETransactions2014}} 
\begin{figure}[htbp]
	\centering
	\subfigure[]{\includegraphics[width=7.5cm,height=2.7cm]{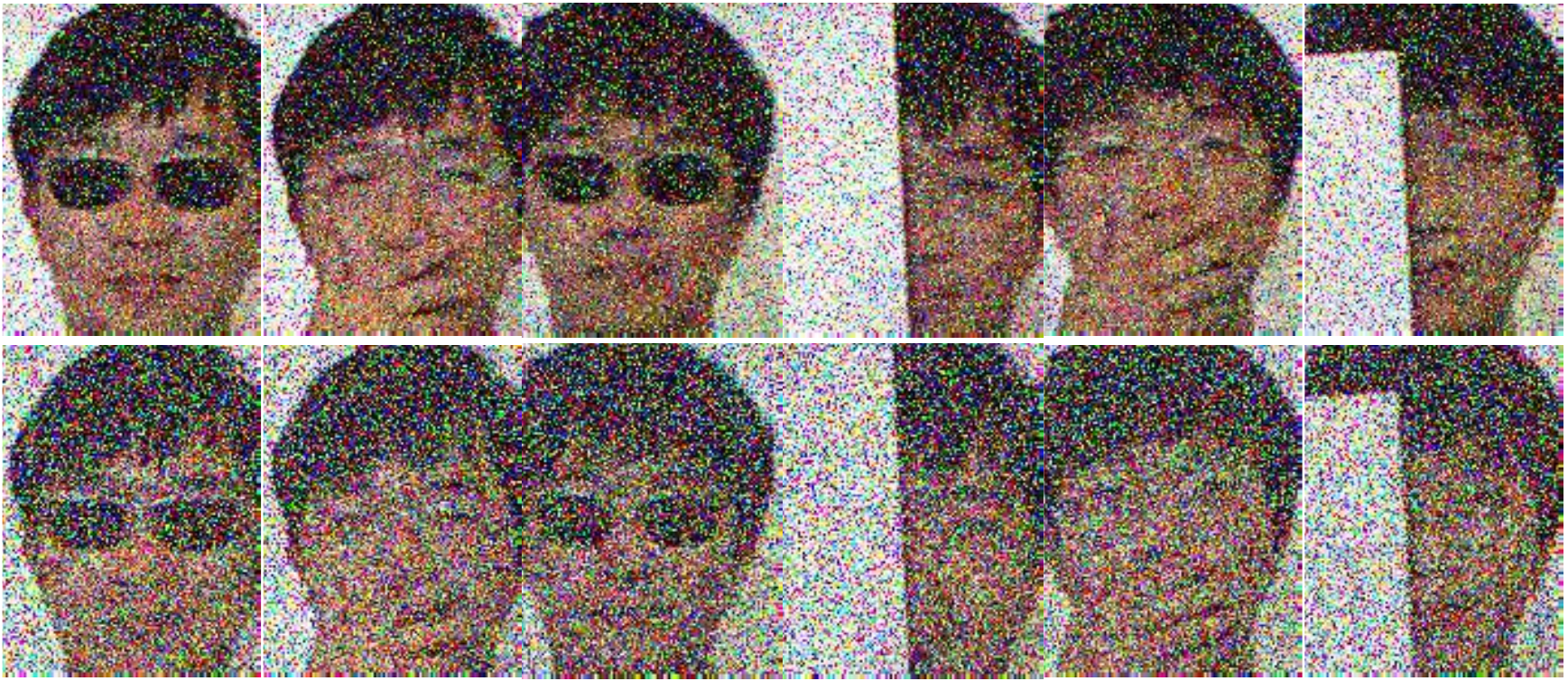}}\\
	\subfigure[]{\includegraphics[width=7.5cm,height=3cm]{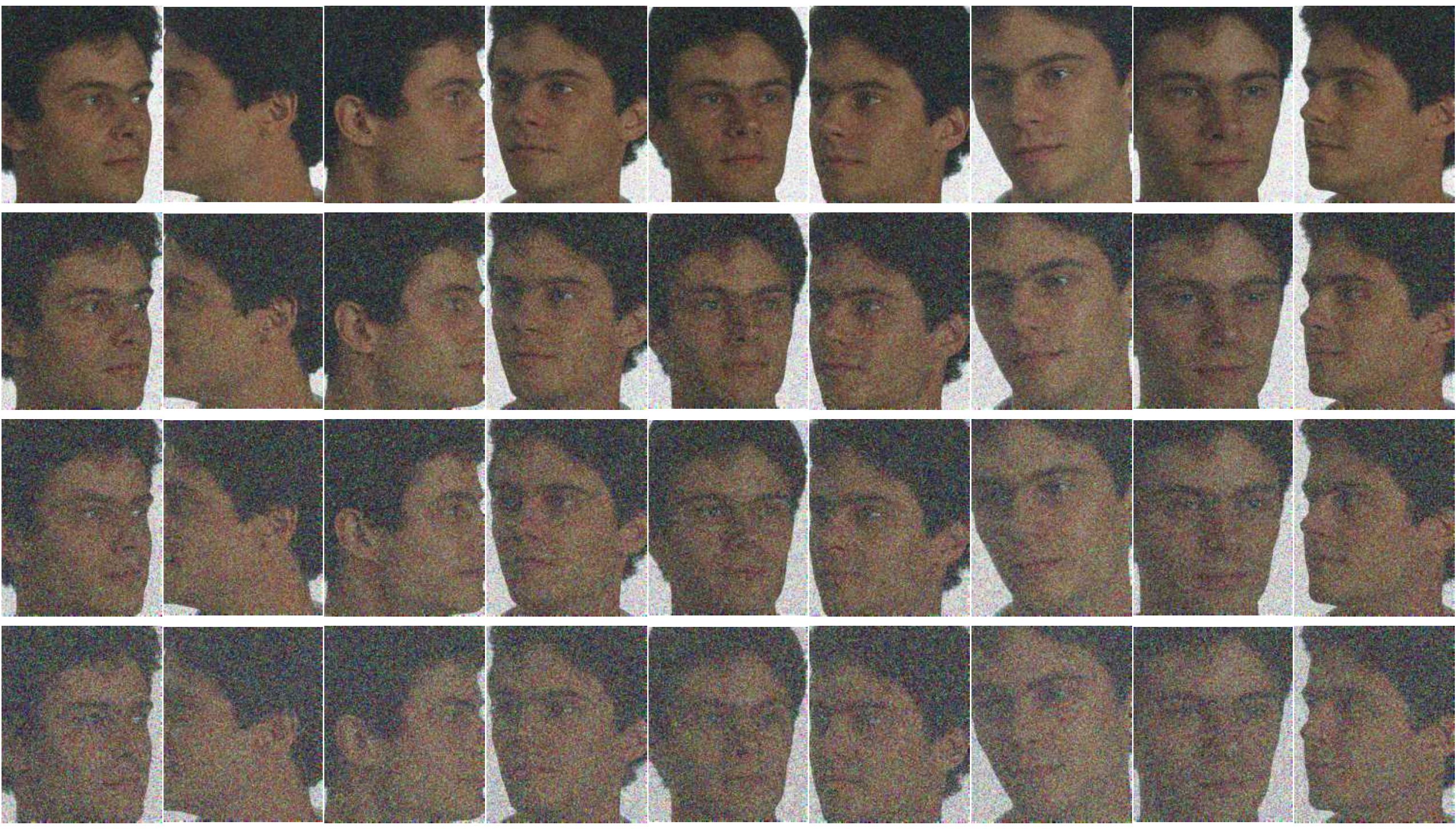}}
	\caption{Sample color face images of one person from (a) EURECOM kinect database (the occluded color face images with different level mixed noise including ``salt \& pepper'' noise and Gaussian noise), (b) Color FERET database face images with different  poses and different level mixed noise including ``salt \& pepper'' noise and Gaussian noise.}
	\label{fig6}
\end{figure}

The EURECOM kinect database consists of the multimodal facial color images of $52$ people ($14$ females, $38$ males) with different facial expressions, different lighting and occlusion conditions. The $14$  frontal face images from both two sessions of each person are used in our experiments. The images of one person from the EURECOM kinect database are shown in Fig. \ref{fig1}(b). All face images are resized to $40\times 40$ pixels. $8$ nonoccluded color face images from two sessions (\emph{e.g.}, the first row in Fig. \ref{fig1}(b)) compose the training set. The test set is formed by adopting three cases. The first case (\textbf{case 1}) 
chooses the $6$ occluded color face images from both two sessions (\emph{e.g.}, the second row in Fig. \ref{fig1}(b)). The second (\textbf{case 2})  and the third (\textbf{case 3})  cases 
still choose the $6$ occluded color face images but with different level mixed noise (Case 2 includes ``salt \& pepper'' noise with $0.2$ probabilities and Gaussian white noise with zero mean and $0.02$ variance.  Case 3 includes ``salt \& pepper'' noise with $0.4$ probabilities and Gaussian white noise with zero mean and $0.04$ variance) from both two sessions. Example test images of cases $2$ and $3$ of one person are respectively shown as the first row and the second row of Fig. \ref{fig6}(a). Based on the results shown in TABLE \ref{tab4}, one can find that the recognition rate of R-NQMR is higher than that of other methods.

\begin{table*}	
	\caption{Recognition rates (\%) of different methods on the EURECOM kinect database for the three cases (the best results are highlighted in bold for each case).}
	\centering
	\begin{tabular}{|c|c|c|c|c|c|c|c|c|c|c|c|}
		\hline 
		\diagbox{Cases}{Rates (\%)}{Methods}&LRC&  QLRC&  CRC&  QCRC&  CROC& QCROC & NMR & RMR & S-RMR & NQMR & R-NQMR \\ 
		\hline 
		case 1&70.19 &70.83  &72.44  &72.44  &76.28  &77.24 &83.33 &90.71 &90.71 &84.29 &\textbf{92.31}  \\ 
		\hline 
		case 2&66.67&67.31 &70.51&71.47&73.72 &74.68&81.41&88.14&88.78&81.73& \textbf{89.10} \\ 
			\hline 
		case 3&48.08&47.76 &67.95&69.87&61.22 &64.74&73.72&86.21&86.54&78.53& \textbf{87.17} \\ 
		\hline 
	\end{tabular} 
	\label{tab4}
\end{table*}

\subsection{Experiments on the Color FERET database \cite{DBLP:journals/pami/PhillipsMRR00}} 
The Color FERET database contains $14, 126$ color face images of $1, 199$ subjects.
The images of one person from Color FERET are shown in Fig. \ref{fig1}(c). We randomly collect a subset that contains $100$ subjects in our experiments. All face images are resized to $32\times 32$ pixels. Each subject has two frontal images
with the letter code ``fa'' and ``fb'' (\emph{e.g.}, the second row in Fig. \ref{fig1}(c)), which are used to compose the training set. The $9$ different poses images (\emph{e.g.}, the first row in Fig. \ref{fig1}(c)) are used for testing. In addition, the mixed noise is added to the images of the test set with the following four cases:
\begin{enumerate}
	\item [1)]: \textbf{Case $1$} includes ``salt \& pepper'' noise with $0.1$ probabilities and Gaussian white noise with zero mean and $0.01$ variance (\emph{e.g.}, the first row in Fig. \ref{fig6}(b)).
	\item [2)]: \textbf{Case $2$} includes ``salt \& pepper'' noise with $0.2$ probabilities and Gaussian white noise with zero mean and $0.02$ variance (\emph{e.g.}, the second row in Fig. \ref{fig6}(b)).
	\item [3)]: \textbf{Case $3$} includes ``salt \& pepper'' noise with $0.3$ probabilities and Gaussian white noise with zero mean and $0.03$ variance (\emph{e.g.}, the third row in Fig. \ref{fig6}(b)).
	\item [4)]: \textbf{Case $4$} includes ``salt \& pepper'' noise with $0.4$ probabilities and Gaussian white noise with zero mean and $0.04$ variance (\emph{e.g.}, the fourth row in Fig. \ref{fig6}(b)).
\end{enumerate}

\begin{table*}	
	\caption{Recognition rates (\%) of different methods on the Color FERET database for the four cases (the best results are highlighted in bold for each case).}
	\centering
	\begin{tabular}{|c|c|c|c|c|c|c|c|c|c|c|c|}
		\hline 
		\diagbox{Cases}{Rates (\%)}{Methods}&LRC&  QLRC&  CRC&  QCRC&  CROC& QCROC & NMR & RMR & S-RMR & NQMR & R-NQMR \\ 
		\hline 
		case 1&55.56 &56.67  &53.78  &54.00  &59.33  &59.11 &58.67 &58.89 &59.11 &59.11 &\textbf{62.22}  \\ 
		\hline 
		case 2&52.22&53.33 &52.22&53.33&54.22 &56.22&57.56&57.33&58.00&58.67& \textbf{62.00} \\ 
		\hline 
		case 3&44.67&47.11 &50.67&50.44&48.89&50.67&56.44&56.89&57.11&57.11& \textbf{61.11} \\ 
				\hline 
		case 4&30.22&36.22&50.00&50.22&37.56&42.67&55.78&53.11&54.22&55.33& \textbf{58.89} \\ 
		\hline 
	\end{tabular} 
	\label{tab5}
\end{table*}
Since the test face images have different poses, all the regression-based algorithms seem to have comparatively low recognition rates shown as TABLE \ref{tab5}.  Even though,  R-NQMR has better performance than that of other methods among all the cases.
\section{Conclusions}
\label{sec6}
Focusing on color face recognition problems, this paper utilizing quaternions to represent the color pixels with RGB channels proposed a nuclear norm based quaternion matrix regression (NQMR) method, which can characterize the low rank structure of the reconstruction error image. Furthermore, to more effectively approximate the matrix rank, a robust nuclear norm based quaternion matrix regression (R-NQMR) method is proposed by adaptively assigning weights on different singular values. Another merit  of R-NQMR method is that it can effectively alleviate the impact of mixed noise in the image on recognition. The ADMM method known that it can guarantee the convergence of the algorithm is developed for solving the two models. Experiments on color face recognition verified that R-NQMR is more effective than some state-of-the-art regression based methods especially for occlusions and noise circumstances. Nevertheless, as with most existing regression based methods, NQMR and R-NQMR still do not recognize faces with large pose changes well, which need further investigation in our future work.

\appendices
\section{Proof of the Theorem \ref{theorem1}}
\label{appendices1}
The proofs of $1)-4)$ can be found in \cite{DBLP:journals/tip/ZouKW16}. Next, we just prove $5)$ and $6)$.

Proof: $5)$
Since,
\begin{align*}\small
\label{app1}
\|\dot{\mathbf{Q}}\|_{F}\!\!=&\left(\sum_{m=1}^{M}\sum_{n=1}^{N}|\dot{q}_{mn}|^{2}\right)^{\frac{1}{2}}\nonumber\\
=&\!\left(\!\sum_{m=1}^{M}\sum_{n=1}^{N}\!|{q_{mn}}_{0}\!+\!{q_{mn}}_{1}i
\!+\!{q_{mn}}_{2}j\!+\!{q_{mn}}_{3}k|^{2}\!\right)^{\frac{1}{2}}\nonumber\\
=&\left(\sum_{m=1}^{M}\sum_{n=1}^{N}\!|{q_{mn}^{2}}_{0}+{q_{mn}^{2}}_{1}
+{q_{mn}^{2}}_{2}+{q_{mn}^{2}}_{3}|^{2}\!\right)^{\frac{1}{2}}\nonumber\\
=&\left(\sum_{m=1}^{M}\sum_{n=1}^{N}\!|{q_{mn}^{2}}_{0}|\right)^{\frac{1}{2}}+\left(\sum_{m=1}^{M}\sum_{n=1}^{N}\!|{q_{mn}^{2}}_{1}|\right)^{\frac{1}{2}}\nonumber\\
&+\left(\sum_{m=1}^{M}\sum_{n=1}^{N}\!|{q_{mn}^{2}}_{2}|\right)^{\frac{1}{2}}+\left(\sum_{m=1}^{M}\sum_{n=1}^{N}\!|{q_{mn}^{2}}_{3}|\right)^{\frac{1}{2}}\nonumber\\
=&\|\mathbf{Q}_{0}\|_{F}+\|\mathbf{Q}_{1}\|_{F}+\|\mathbf{Q}_{2}\|_{F}+\|\mathbf{Q}_{3}\|_{F},
\end{align*}
and from (\ref{qtorm}), we have $\|\dot{\mathbf{Q}}\|_{F}=\frac{1}{4}\|\mathcal{P}(\dot{\mathbf{Q}})\|_{F}$.

$6)$ We first give the following Lemma \ref{lemma2}:
\begin{lemma}
\label{lemma2}	
Let $\dot{\mathbf{P}}$ and $\dot{\mathbf{Q}}\in\mathbb{H}^{M\times N}$, then
\begin{enumerate}
	\item [\textcircled{1}] $\mathcal{P}(\dot{\mathbf{P}}\dot{\mathbf{Q}})=\mathcal{P}(\dot{\mathbf{P}})\mathcal{P}(\dot{\mathbf{Q}})$.
	\item [\textcircled{2}] $\mathcal{P}(\dot{\mathbf{Q}})^{T}=\mathcal{P}(\dot{\mathbf{Q}}^{H})$.
	\item [\textcircled{3}] If $\dot{\mathbf{P}}$ has unitary columns (\emph{i.e.}, $\dot{\mathbf{P}}^{H}\dot{\mathbf{P}}=\mathbf{I}$ \cite{10029950538}), then $\mathcal{P}(\dot{\mathbf{P}})$ has orthogonal columns.
\end{enumerate}
\end{lemma}

Proof: After some trivial matrix calculations, the $\textcircled{1}$ is obviously true. For $\textcircled{2}$, we see that 
\begin{align*}\small
\mathcal{P}(\dot{\mathbf{Q}})^{T}&=\left(\begin{array}{cccc}
\mathbf{Q}_{0}^{T}&  \mathbf{Q}_{1}^{T}&  \mathbf{Q}_{2}^{T} &\mathbf{Q}_{3}^{T} \\ 
\mathbf-{Q}_{1}^{T}&  \mathbf{Q}_{0}^{T}&  \mathbf{Q}_{3}^{T} &-\mathbf{Q}_{2}^{T} \\ 
-\mathbf{Q}_{2}^{T}&  -\mathbf{Q}_{3}^{T}& \mathbf{Q}_{0}^{T} &\mathbf{Q}_{1}^{T}  \\ 
-\mathbf{Q}_{3}^{T}&  \mathbf{Q}_{2}^{T}&  -\mathbf{Q}_{1}^{T} &\mathbf{Q}_{0}^{T}
\end{array} \right)\\
&=\mathcal{P}\big(\mathbf{Q}_{0}^{T}-\mathbf{Q}_{1}^{T}i-\mathbf{Q}_{2}^{T}j-\mathbf{Q}_{3}^{T}k\big)\\
&=\mathcal{P}(\dot{\mathbf{Q}}^{H}).
\end{align*}
For $\textcircled{3}$, based on $\textcircled{1}$ and $\textcircled{2}$, we have $\mathcal{P}(\dot{\mathbf{Q}})^{T}\mathcal{P}(\dot{\mathbf{Q}})=\mathcal{P}(\dot{\mathbf{Q}}^{H})\mathcal{P}(\dot{\mathbf{Q}})=\mathcal{P}(\dot{\mathbf{Q}}^{H}\dot{\mathbf{Q}})=\mathcal{P}(\mathbf{I})=\mathbf{I}$, \emph{i.e.}, $\mathcal{P}(\dot{\mathbf{P}})$ has orthogonal columns.

Assume the SVD of the quaternion matrix $\dot{\mathbf{Q}}$ with rank $r$ is
\begin{equation*}
\dot{\mathbf{Q}}=\dot{\mathbf{U}}\mathbf{\Sigma}_{r}\dot{\mathbf{V}}^{H},
\end{equation*} 
where $\dot{\mathbf{U}}\in\mathbb{H}^{M\times r}$ and $\dot{\mathbf{V}}\in\mathbb{H}^{N\times r}$ have unitary columns, and $\mathbf{\Sigma}_{r}={\rm{diag}}(\sigma_{1},\ldots, \sigma_{r} )\in\mathbb{R}^{r\times r}$ consists of all positive singular values of quaternion matrix $\dot{\mathbf{Q}}$. Then, based on Lemma \ref{lemma2}, we have $\mathcal{P}(\dot{\mathbf{Q}})=\mathcal{P}(\dot{\mathbf{U}})\mathcal{P}(\mathbf{\Sigma}_{r})\mathcal{P}(\dot{\mathbf{V}})^{T}$ which is obvious one of the SVDs of matrix $\mathcal{P}(\dot{\mathbf{Q}})$. And we can find that the $\mathcal{P}(\mathbf{\Sigma}_{r})$ has the following form:
\begin{equation}\small
\label{app2}
\mathcal{P}(\mathbf{\Sigma}_{r})=\left(\begin{array}{cccc}
\mathbf{\Sigma}_{r}&  \mathbf{0}&  \mathbf{0}& \mathbf{0} \\ 
\mathbf{0}&  \mathbf{\Sigma}_{r}&  \mathbf{0} & \mathbf{0} \\ 
\mathbf{0}&  \mathbf{0}&   \mathbf{\Sigma}_{r}& \mathbf{0} \\ 
\mathbf{0}&  \mathbf{0}&  \mathbf{0}& \mathbf{\Sigma}_{r}
\end{array} \right).
\end{equation}
Thus, we have $\|\mathcal{P}(\dot{\mathbf{Q}})\|_{\ast}=4\sum_{k=1}^{r}\sigma_{k}=4\|\dot{\mathbf{Q}}\|_{\ast}$.

\section*{Acknowledgment}

This work was supported by The Science and Technology Development Fund, Macau SAR (File no. FDCT/085/2018/A2).

\ifCLASSOPTIONcaptionsoff
  \newpage
\fi

\bibliographystyle{IEEEtran}
\bibliography{Myreference}
\end{document}